\documentclass[11pt,a4paper]{article}
\usepackage{a4,amssymb,latexsym,amsmath,color}
\begin{document}


\newcounter{pkt}
\newenvironment{enum}{\setcounter{pkt}{0}
\begin{list}{\rm\alph{pkt})}{\usecounter{pkt}
\setlength{\topsep}{1ex}\setlength{\labelwidth}{0.5cm}
\setlength{\leftmargin}{1cm}\setlength{\labelsep}{0.25cm}
\setlength{\parsep}{-3pt}}}{\end{list}~\\[-6ex]}

\newcounter{punkt1}
\newenvironment{enum1}{\setcounter{punkt1}{0}
\begin{list}{(\arabic{punkt1})}{\usecounter{punkt1}
\setlength{\topsep}{1ex}\setlength{\labelwidth}{0.6cm}
\setlength{\leftmargin}{1cm}\setlength{\labelsep}{0.25cm}
\setlength{\parsep}{-3pt}}}{\end{list}~\\[-6ex]}

\newcounter{punkt2}
\newenvironment{enum2}{\setcounter{punkt2}{0}
\begin{list}{(\roman{punkt2})}{\usecounter{punkt2}
\setlength{\topsep}{1ex}\setlength{\labelwidth}{0.6cm}
\setlength{\leftmargin}{1cm}\setlength{\labelsep}{0.25cm}
\setlength{\parsep}{-3pt}}}{\end{list}~\\[-6ex]}

\newcounter{punkt3}
\newenvironment{enumbib}{\setcounter{punkt3}{0}
\begin{list}{\arabic{punkt3}.}{\usecounter{punkt3}
\setlength{\topsep}{1ex}\setlength{\labelwidth}{0.6cm}
\setlength{\leftmargin}{1cm}\setlength{\labelsep}{0.25cm}
\setlength{\parsep}{1pt}}}{\end{list}}

\newenvironment{bew}
{\vspace*{-0.25cm}\begin{sloppypar}\noindent{\it 
Proof.}}{\hfill\qed\end{sloppypar}\vspace*{0.15cm}}

\newcounter{pic}\setcounter{pic}{0}

\newcommand{\bdpm}{\begin{displaymath}}
\newcommand{\edpm}{\end{displaymath}}

\newcommand{\beas}{\begin{eqnarray*}}
\newcommand{\eeas}{\end{eqnarray*}}

\newcommand{\ba}{\begin{array}}
\newcommand{\ea}{\end{array}}

\newtheorem{theo}{Theorem}[section]          
\newtheorem{lem}[theo]{Lemma}
\newtheorem{cor}[theo]{Corollary}
\newtheorem{prop}[theo]{Proposition}
\newtheorem{exam}[theo]{Example}
\newtheorem{rem}[theo]{Remark}
\newtheorem{rems}[theo]{Remarks}

\newcommand{\brm}{\begin{rm}}
\newcommand{\erm}{\end{rm}}

\newcommand{\qed}{\hfill $\Box$}
\newcommand{\des}{{\sf des}}
\newcommand{\ddes}{{\sf ddes}}
\newcommand{\exc}{{\sf exc}}
\newcommand{\dexc}{{\sf dexc}}
\newcommand{\inv}{{\sf inv}}
\newcommand{\den}{{\sf den}}
\newcommand{\maj}{{\sf maj}}
\newcommand{\rk}{{\rm rank}}
\newcommand{\nat}{{\Bbb N}}


\begin{center}
{\large\bf THE EXCEDANCES AND DESCENTS}\\ 
{\large\bf OF BI-INCREASING PERMUTATIONS}\\[1cm]
Astrid Reifegerste\\
Institut f\"ur Mathematik, Universit\"at Hannover\\
Welfengarten 1, D-30167 Hannover, Germany\\
{\it reifegerste@math.uni-hannover.de}\\[0.5cm]
December 16, 2002
\end{center}
\vspace*{0.5cm}

\begin{footnotesize}
{\sc Abstract.} 
Starting from some considerations we make about the relations between certain 
difference statistics and the classical permutation statistics we study permutations whose  
inversion number and excedance difference coincide. It turns out that 
these (so-called bi-increasing) permutations are just the $321$-avoiding ones.
The paper investigates their excedance and descent structure. In particular, we 
find some nice combinatorial interpretations for the distribution coefficients of the 
number of excedances and descents, respectively, and their difference analogues over the 
bi-increasing permutations in terms of parallelogram polyominoes and 2-Motzkin 
paths. This yields a connection between restricted permutations, parallelogram 
polyominoes, and lattice paths that reveals the relations between several 
well-known bijections given for these objects (e.g. by Delest-Viennot, 
Billey-Jockusch-Stanley, Fran\c{c}on-Viennot, and 
Foata-Zeilberger). As an application, we enumerate skew diagrams according to their 
rank and give a simple combinatorial proof for a result concerning the symmetry 
of the joint distribution of the number of excedances and inversions, respectively, over the 
symmetric group.
\end{footnotesize}
\vspace*{1cm}


\setcounter{section}{1}\setcounter{theo}{0}

\centerline{\large{\bf 1}\hspace*{0.25cm}
{\sc Motivation and preliminaries}}
\vspace*{0.5cm}

Let ${\cal S}_n$ be the set of all permutations of $[n]:=\{1,\ldots,n\}$. We 
write any permutation $\pi\in{\cal S}_n$ as word $\pi_1\cdots\pi_n$ where 
$\pi_i$ means the integer $\pi(i)$.\\[2ex]
In this paper, we will mainly investigate certain permutations in view of the 
behaviour of their excedances and descents. First we recall the definitions and 
fix some notations.\\[2ex]  
For $\pi\in{\cal S}_n$, an {\it excedance} of $\pi$ is an integer $i\in[n-1]$ such that $\pi_i>i$. Here 
the element $\pi_i$ is called an {\it excedance letter}. The set of excedances 
of $\pi$ is denoted by ${\sf E}(\pi)$. By $\pi_{\sf e}$ and $\pi_{\sf ne}$ we denote the restrictions 
of $\pi$ on the excedances and non-excedances, respectively.\\
A {\it descent} of $\pi$ is an integer $i\in[n-1]$ for which $\pi_i>\pi_{i+1}$. 
If $i$ is a descent we call $\pi_i$ a {\it descent top} and $\pi_{i+1}$ a {\it descent 
bottom}. The set of descents of $\pi$ is denoted by ${\sf D}(\pi)$. We write $\pi_{\sf d}$ to denote the
subword consisting of all descent tops of $\pi$ (in order of appearance). The 
subword formed from the remaining letters we denote by $\pi_{\sf nd}$.\\
A pair $(\pi_i,\pi_j)$ is called an {\it inversion} of $\pi$ if $i<j$ and 
$\pi_i>\pi_j$. The set containing all inversions of $\pi$ we denote by ${\sf 
I}(\pi)$. (Clearly, inversions are defined in the same way for arbitrary 
words.)\\
Three of the four classical statistics count the number of occurrences of these 
patterns in a permutation, the fourth one, the so-called {\it major index}, is the sum of descents. We use the 
usual notation: 
\bdpm
\exc(\pi)=|{\sf E}(\pi)|,\quad\des(\pi)=|{\sf D}(\pi)|,\quad
\inv(\pi)=|{\sf I}(\pi)|,\quad\maj(\pi)=\sum_{i\in{\sf D}(\pi)} i.
\edpm
We write $S_n^{\sf f}(k)$ to denote the number of permutations in ${\cal S}_n$ 
for which the statistic {\sf f} takes the value $k$. Analogously, 
the coefficients of the joint distribution of {\sf f} and {\sf g} on ${\cal 
S}_n$ are denoted by $S_n^{\sf (f,g)}(k,l)$. When we consider statistics over 
the set ${\cal B}_n$ defined below then we write $B$ instead of $S$.\\
In \cite{clarke etal}, the authors studied some differences of permutation statistics. 
Here we will deal with two of these, namely the {\it excedance 
difference} and the {\it descent difference} of a permutation $\pi\in{\cal S}_n$ defined by 
\bdpm
\dexc(\pi)=\sum_{i\in{\sf E}(\pi)} (\pi_i-i) \qquad\mbox{and}\qquad
\ddes(\pi)=\sum_{i\in{\sf D}(\pi)} (\pi_i-\pi_{i+1}),
\edpm
respectively.\\
For example, the permutation $\pi=4\:2\:8\:3\:6\:9\:7\:5\:1\:10\in{\cal S}_{10}$ 
has the excedances  $1,3,5,6$ and the descents $1,3,6,7,8$. Its excedance 
letter word and non-excedance letter word, respectively, are $\pi_{\sf 
e}=4\:8\:6\:9$ and $\pi_{\sf ne}=2\:3\:7\:5\:1\:10$. For the descent top word 
we obtain $\pi_{\sf d}=4\:8\:9\:7\:5$. The excedance difference equals 
\bdpm
\dexc(\pi)=(4-1)+(8-3)+(6-5)+(9-6)=12,
\edpm
the descent difference is equal to 
\bdpm
\ddes(\pi)=(4-2)+(8-3)+(9-7)+(7-5)+(5-1)=15.
\edpm
This section is to show that the difference statistics are closely connected 
with the classical ones.\\[2ex]
Foata showed, combinatorially, that the excedance number is {\it 
equidistributed} with the descent number over ${\cal S}_n$. That is to say, 
$S_n^{\exc}(k)=S_n^{\des}(k)$, for all $n$ and $k$. The same result holds for 
the difference statistics. Foata's bijection \cite[Th. 10.2.3]{foata} 
also proves this fact, and hence even the equidistribution of the pairs 
$(\exc,\dexc)$ and $(\des,\ddes)$. For proof we review Foata's proof. 

\begin{prop} \label{foata}
Both statistics $\dexc$ and $\ddes$ and bistatistics $(\exc,\dexc)$ and 
$(\des,\ddes)$ are equidistributed over ${\cal S}_n$.
\end{prop}

\begin{bew}
Given a permutation $\pi\in{\cal S}_n$, decompose $\pi$ in distinct cycles 
\bdpm
\pi=(c_{11}c_{12}\cdots c_{1l_1})(c_{21}c_{22}\cdots c_{2l_2})\cdots
(c_{s1}c_{s2}\cdots c_{sl_s})
\edpm
such that (i) each cycle is written with its largest element first, and (ii) 
the cycles are written in increasing order of their first element. Fixed 
points are regarded as singleton cycle. Define the map 
\bdpm
\varphi:\:{\cal S}_n\to{\cal S}_n,\;\pi\mapsto
c_{11}c_{1l_1}c_{1{l_1-1}}\cdots c_{12}\,c_{21}c_{2l_2}c_{2{l_2-1}}\cdots 
c_{22}\cdots c_{s1}c_{sl_s}c_{s{l_s-1}}\cdots c_{s2}.
\edpm
By Foata, $\phi$ is bijectiv, and we have $\exc(\pi)=\des(\phi(\pi))$.\\
But we have $\dexc(\pi)=\ddes(\phi(\pi))$ as well. To see this, first let $i$ be an excedance of 
$\pi$. Clearly, $i$ and $\pi_i$ belong to the same cycle $c$, and appear in 
$\phi(\pi)$ as consecutive letters, beginning with $\pi_i$. (If $\pi_i$ is the 
largest element of $c$ then $c=(\pi_i\cdots i)$, otherwise $c=(\pi_k\cdots 
i\:\pi_i\cdots)$ for some $k\not=i$.) If $i$ is a fixed point of $\pi$ condition (ii) 
causes that the letter following $i$ in $\phi(\pi)$ is greater than $i$. 
Finally, let $i$ be an integer satisfying $i>\pi_i$. Then the cycle $c$ containing $i$ 
and $\pi_i$ is written as either $c=(k\cdots i\:\pi_i\cdots)$ for some $k\not=i$, and 
$\pi_i\,i$ occurs in $\phi(\pi)$, or $c=(i\:\pi_i\:\pi_i^2\cdots\pi_i^{l-1})$, 
and $i\,\pi_i^{l-1}\pi_i^{l-2}\cdots\pi_i$ occurs in $\phi(\pi)$ where $\pi_i^k$ denotes $\pi(\pi_i^{k-1})$ for $k\ge 
2$. (Note that $i$ is the largest cycle element in this case.)\\
Consequently, the positive differences $\pi_i-i$ correspond exactly to the 
positive differences between consecutive letters in $\phi(\pi)$.   
\end{bew}
\vspace*{2ex}

Because of the definition, the relation between the difference statistics and 
their classical counterparts comes as no surprise; more amazing is the fact that 
$\dexc$ is closely connected with the number of inversions, as well. This 
result already appeared in \cite[Th. 2]{clarke etal}.\\
We give another but just as simple proof. Its key object is the permutation $\hat{\pi}\in{\cal S}_n$ 
that is obtained from $\pi$ by sorting the letters of $\pi_{\sf e}$ 
and $\pi_{\sf ne}$ in increasing order, respectively. (That is to mean, 
$\hat{\pi}_{i_k}=\pi_{\sigma(i_k)}$ and $\hat{\pi}_{j_k}=\pi_{\tau(j_k)}$ where 
$i_1,\ldots,i_e$ are the excedances, and $j_1,\ldots,j_{n-e}$ the 
remaining elements of $\pi$, and $\sigma$ and $\tau$ are permutations of ${\sf 
E}(\pi)$ and the set of non-excedances, respectively, such that 
$\hat{\pi}_{i_1}<\ldots<\hat{\pi}_{i_e}$ and 
$\hat{\pi}_{j_1}<\ldots<\hat{\pi}_{j_{n-e}}$.)
For example, let $\pi=\underline{4}\:2\:\underline{8}\:3\:\underline{6}\:\underline{9}\:7\:5\:1\:10\in{\cal S}_{10}$ as before, with $\pi_{\sf 
e}=4\:8\:6\:9$ and $\pi_{\sf ne}=2\:3\:7\:5\:1\:10$. (The letters of the 
subword $\pi_{\sf e}$ are underlined.) Then we obtain the bi-sorted permutation 
$\hat{\pi}=\underline{4}\:1\:\underline{6}\:2\:\underline{8}\:\underline{9}\:3\:5\:7\:10$.

\begin{theo} \label{relation inv - dexc}
For any permutation $\pi\in{\cal S}_n$ we have $\inv(\pi)=\dexc(\pi)+\inv(\pi_{\sf e})+\inv(\pi_{\sf 
ne})$.
\end{theo}

\begin{bew}
First we show that $\exc(\hat{\pi})=\exc(\pi)$ and $\dexc(\hat{\pi})=\dexc(\pi)$.\\
Any excedance of $\pi$ is an excedance in $\hat{\pi}$ as well: let $i,j\in{\sf 
E}(\pi)$ such that $i<j$ and $\pi_i>\pi_j$. Obviously, $i<j<\pi_j$ and 
$j<\pi_j<\pi_i$. Analogously, any non-excedance of $\pi$ is a non-excedance in 
$\hat{\pi}$. Thus the permutations $\pi$ and $\hat{\pi}$ have the same ecxedance 
set, and we obtain
\bdpm
\dexc(\pi)=\sum_{i\in{\sf E}(\pi)} (\pi_i-i)=\sum_{i\in{\sf E}(\pi)} 
(\pi_{\sigma(i)}-i)=\sum_{i\in{\sf E}(\hat{\pi})} (\hat{\pi}_i-i)=\dexc(\hat{\pi}).
\edpm
Now we determine the inversion number of $\hat{\pi}$. To this end, we utilize the {\it 
reverse code} $c(\pi)=(c_1,\ldots,c_n)$ of $\pi$ whose $i$th component is defined 
to be the number of integers $j\in[n]$ satisfying $j<i$ and $\pi_j>\pi_i$. (Clearly, 
$c_1+\ldots+c_n=\inv(\pi)$.)\\
Let $i_1,\ldots,i_e$ be the excedances of $\pi$ where 
$\pi_{i_1}<\pi_{i_2}<\ldots<\pi_{i_k}>\pi_{i_{k+1}}$. Then $c_{i_k}=0$ and 
$c_{i_{k+1}}\ge1$. For the reverse code $(c'_1,\ldots,c'_n)$ of the permutation obtained by 
exchanging $\pi_{i_k}$ and $\pi_{i_{k+1}}$ 
we have 
\beas
&&c'_{i_k}=c_{i_{k+1}}-1,\\
&&c'_{i_{k+1}}=0,\\
&&c'_j=c_j\quad\mbox{for all }j\in[n]\setminus\{i_k,i_{k+1}\}.
\eeas
The last relation is evident for 
$j=1,\ldots,i_k-1,i_{k+1}+1,\ldots,n$. By definition, $i_k+1,\ldots,i_{k+1}-1$ 
are non-excedances each of which satisfies $\pi_j\le 
j<i_{k+1}<\pi_{i_{k+1}}<\pi_{i_k}$. Thus the inversion $(\pi_{i_k},\pi_j)$ 
appearing in $\pi$ corresponds to $(\pi_{i_{k+1}},\pi_j)$ in $\hat{\pi}$. Consequently, 
exchanging the excedance letters $\pi_{i_k}$ and $\pi_{i_{k+1}}$ decreases the number of inversions by 1. 
(Sorting $\pi_{\sf e}$ completely requires $\inv(\pi_{\sf e})$ transpositions.) By 
reasoning similar for non-excedances, we obtain 
\bdpm
\inv(\hat{\pi})=\inv(\pi)-\inv(\pi_{\sf e})-\inv(\pi_{\sf ne}).
\edpm
By construction, the words $\hat{\pi}_{\sf e}$ and $\hat{\pi}_{\sf 
ne}$ do not contain any inversion. Hence the component $\hat{c}_i$ of 
$\hat{\pi}$'s reverse code equals zero if $i$ is an excedance of $\hat{\pi}$ (or 
equivalently, of $\pi$), and is equal to the difference $i-\hat{\pi}_i$ 
otherwise. Thus 
\bdpm
\inv(\hat{\pi})=\sum_{i=1}^n \hat{c}_i=\sum_{i\notin{\sf E(\hat{\pi})}} 
(i-\hat{\pi}_i)=\sum_{i\in{\sf E(\hat{\pi})}} (\hat{\pi}_i-i)=\dexc(\hat{\pi}),
\edpm
and the assertion is proved.
\end{bew}
\vspace*{2ex}

Define two permutations $\pi,\sigma\in{\cal S}_n$ as being {\it equivalent} if 
$\hat{\pi}=\hat{\sigma}$. The fundamental property that both $\exc$ and $\dexc$ 
are constant on the equivalence classes we will use in Section 6 to derive 
results concerning the distribution of these statistics over ${\cal S}_n$ from 
analogous results proved for the set of class representatives.\\[2ex]
The theorem yields an estimate for the inversion number depending 
on the statistics $\exc$ and $\dexc$. Its proof needs the first one of the following auxiliary 
results.

\begin{lem} \label{auxiliary}
\begin{enum}
\item[]
\item Let $\pi\in{\cal S}_n$ satisfy $\inv(\pi)=\dexc(\pi)$, and let $k$ be a non-excedance 
of $\pi$. Then 
\bdpm
|\{j\notin{\sf E}(\pi):\pi_k\le j<k\}|=|\{i\in{\sf E}(\pi):i<\pi_k<\pi_i\}|.
\edpm
\item Let $\pi\in{\cal S}_n$ be an arbitrary permutation, and let $k$ be a non-excedance 
of $\pi$. Then 
\bdpm
|\{j\notin{\sf E}(\pi):\pi_j\le k<j\}|=|\{i\in{\sf E}(\pi):i<k<\pi_i\}|.
\edpm
\end{enum}
\end{lem}
\vspace*{-3ex}

\begin{bew}
{\bf a)} By the theorem, both $\pi_{\sf e}$ and $\pi_{\sf ne}$ are increasing 
words. Thus the 
letters $1,2,\ldots,\pi_k-1$ appear to the left of position $k$. Consequently, 
$k-\pi_k$ integers $i<k$ satisfy $\pi_i>\pi_k$ (here $i$ is necessarily an excedance of $\pi$). 
Let $s$ be the number of these integers for which $i<\pi_k$ in addition. Hence 
among the integers $\pi_k,\pi_k+1,\ldots,k-1$ there are $k-\pi_k-s$ excedances, 
and therefore, $s$ non-excedances.\\
{\bf b)} Evidently, $k-1$ letters $\pi_j$ different from $\pi_k$ satisfy 
$\pi_j\le k$. Let $t$ be the number of all these positions $j$ being to the 
right of $k$. Clearly, such a integer $j$ is a non-excedance of $\pi$. In 
$[k-1]$ there are $k-1-t$ integers $i$ for which $\pi_i\le k$, and 
hence $t$ ones with $\pi_i>k$. Obviously, the latter ones are necessarily excedances.
\end{bew}

\begin{cor} \label{bounds}
For any permutation $\pi\in{\cal S}_n$ we have $\dexc(\pi)\le\inv(\pi)\le2\dexc(\pi)-\exc(\pi)$.
\end{cor}

\begin{bew}
The lower bound follows immediately from Theorem \ref{relation inv - dexc}. To 
prove the upper one we first construct a permutation sequence 
$\pi^{(0)},\pi^{(1)},\ldots$ satisfying $\dexc(\pi^{(k)})=k$ and 
$\inv(\pi^{(k)})=2k-\exc(\pi^{(k)})$. Initialize $\pi^{(0)}=12\cdots n$. Given 
$\pi^{(k-1)}$, we define $\pi^{(k)}$ as follows:
\begin{enum1}
\item Let $e$ be the smallest integer such that $e+\pi_e^{(k-1)}<n+1$. If there 
is no such an $e$ stop the procedure. Set $a:=\pi_e^{(k-1)}$.
\item Set $\pi_e^{(k)}:=a+1$, $\pi_a^{(k)}:=a$ if $a\not=e$, 
$\pi_{a+1}^{(k)}:=e$, and $\pi_i^{(k)}:=\pi_i^{(k-1)}$ otherwise.
\end{enum1}
\vspace*{-0.3cm}

(For an example, see the following remark.) Note that $\exc(\pi^{(k)})=e$; more 
exactly, we have ${\sf E}(\pi^{(k)})=[e]$. By the construction, 
$\dexc$ increases with $k$. The inversion number is 
increased by 2 if the integer $e$ does not change with regard to the previous 
step. Otherwise we have $\inv(\pi^{(k)})=\inv(\pi^{(k-1)})+1$. Hence 
$\inv(\pi^{(k)})=2k-e$.\\
On the other hand, this is also the maximum inversion number for a permutation $\pi\in{\cal S}_n$ with $e$ excedances and 
excedance difference $k$. Why? By Theorem \ref{relation inv - dexc}, it suffices to show $\inv(\pi_{\sf e})+\inv(\pi_{\sf 
ne})\le k-e$. As shown in the proof of \ref{relation inv - dexc}, this assertion is equivalent to the following 
one. Let $\sigma\in{\cal S}_n$ be a permutation for which $\exc(\pi)=e$ and $\inv(\pi)=\dexc(\pi)=k$. 
Then there are at most $k-e$ inversions on the excedance letters and 
non-excedance letters of $\sigma$ which preserve the set of excedances.\\ 
Let $k<i$ be excedances of $\sigma$. (Note that 
$\sigma_k<\sigma_i$.) Then the exchange of $\sigma_k$ and $\sigma_i$ keeps all 
excedances if and only if $i<\sigma_k$. Therefore, for any $\pi\in{\cal S}_n$ satisfying 
$\hat{\pi}=\sigma$ we have
\bdpm
\inv(\pi_{\sf e})\le\sum_{k\in{\sf E}(\sigma)} |\{i\in{\sf 
E}(\sigma):i<\sigma_k<\sigma_i\}|.
\edpm 
Now let $j<k$ be non-excedances of $\sigma$. (We have $\sigma_j<\sigma_k$ 
again.) When exchanging $\sigma_j$ and $\sigma_k$ the integers $j$ and $k$ are 
still non-excedances if and only if $j\ge\sigma_k$. By Lemma \ref{auxiliary}a, 
the number of non-excedances $j$ satisfying $\sigma_k\le j<k$ equals the number 
of excedances $i$ for which $i<\sigma_k<\sigma_i$. Consequently, we have 
\bdpm
\inv(\pi_{\sf ne})\le\sum_{k\notin{\sf E}(\sigma)} |\{i\in{\sf 
E}(\sigma):i<\sigma_k<\sigma_i\}|,
\edpm 
for each $\pi\in{\cal S}_n$ with $\hat{\pi}=\sigma$. Both together yield
\bdpm
\inv(\pi_{\sf e})+\inv(\pi_{\sf ne})\le\sum_k |\{i\in{\sf 
E}(\sigma):i<\sigma_k<\sigma_i\}|=\sum_{i\in{\sf E}(\sigma)} 
(\sigma_i-i-1)=\dexc(\sigma)-\exc(\sigma)
\edpm
which is the desired bound.
\end{bew}

\begin{cor}
The conditions $\exc(\pi)=\inv(\pi)$ and $\exc(\pi)=\dexc(\pi)$ are equivalent 
in ${\cal S}_n$.
\end{cor}

\begin{rems}
\brm
\begin{enum}
\item[]
\item 
For $n=6$, the procedure given in the proof of \ref{bounds} yields the sequence:
\vspace*{1ex}

\begin{footnotesize}
\begin{center}
\begin{tabular}{|c|c|c|c|c|}\hline
$k$&$\pi^{(k)}$&$\exc(\pi^{(k)})$&$\dexc(\pi^{(k)})$&$\inv(\pi^{(k)})$\\\hline
0&$1\:2\:3\:4\:5\:6$&0&0&0\\
1&$2\:1\:3\:4\:5\:6$&1&1&1\\
2&$3\:2\:1\:4\:5\:6$&1&2&3\\
3&$4\:2\:3\:1\:5\:6$&1&3&5\\
4&$5\:2\:3\:4\:1\:6$&1&4&7\\
5&$6\:2\:3\:4\:5\:1$&1&5&9\\
6&$6\:3\:2\:4\:5\:1$&2&6&10\\
7&$6\:4\:3\:2\:5\:1$&2&7&12\\
8&$6\:5\:3\:4\:2\:1$&2&8&14\\
9&$6\:5\:4\:3\:2\:1$&3&9&15\\\hline
\end{tabular}  
\end{center}
\end{footnotesize}
\vspace*{1ex}
 
It is easy to see that $\dexc(\pi)\le\frac{1}{4}n^2$ for each permutation $\pi\in{\cal S}_n$. 
How the construction shows, for every integer $0\le k\le\frac{1}{4}n^2$ there is 
a permutation $\pi\in{\cal S}_n$ having excedance difference $k$.
\item Let $I=\{i_1,\ldots,i_e\}$ and $A=\{a_1,\ldots,a_e\}$ be ordered integer 
sets. The proof gives implicit a simple construction for a permutation in ${\cal S}_n$ whose 
excedance set is $I$, whose excedance letter set is $A$, and 
whose inversion number is maximal. Let $A_{i_k}$ be the set consisting 
of all $a\in A$ satisfying $a>i_k$ where $k=1,\ldots,e$. Define $\pi_{i_e}$ to 
be the smallest element of $A_{i_e}$, and delete it from 
$A_{i_1},\ldots,A_{i_{e-1}}$. Then define $\pi_{i_{e-1}}$ to 
be the smallest element of $A_{i_{e-1}}$, and delete it from 
$A_{i_1},\ldots,A_{i_{e-2}}$, and so on. Proceed analogously with the 
non-excedances and their letters. Let $B_{j_k}$ be consisted of the elements 
$b\in B$ satisfying $b\le j_k$ where $B$ is the $[n]$-complement of $A$, and 
$j_1<\ldots<j_{n-e}$ are the elements of $[n]$, different from 
$i_1,\ldots,i_e$. Define $\pi_{j_1}$ to be the largest element of $B_{j_1}$, 
and delete it from $B_{j_2},\ldots,B_{j_{n-e}}$. Then define 
$\pi_{j_2}$ to be the largest element of $B_{j_2}$, and delete it from 
$B_{j_3},\ldots,B_{j_{n-e}}$, and so on.\\ 
For instance, given $I=\{1,3,5,6\}$, $A=\{4,6,8,9\}$, and $n=10$ we obtain   
\beas
&&A_6=\{8,9\},\;A_5=\{6,\not8,9\},\;A_3=\{4,\not6,\not8,9\},\;A_1=\{\not4,\not6,\not8,9\};\\
&&B_2=\{1,2\},\;B_4=\{1,\not2,3\},\;B_7=\{1,\not2,\not3,5,7\},\;B_8=\{1,\not2,\not3,5,\not7\},\\
&&B_9=\{1,\not2,\not3,\not5,\not7\},\;B_{10}=\{\not1,\not2,\not3,\not5,\not7,10\};
\eeas
and hence $\pi=\underline{9}\:2\:\underline{4}\:3\:\underline{6}\:\underline{8}
\:7\:5\:1\:10\in{\cal S}_{10}$ having $20$ inversions.\\
Clearly (by the proof of Theorem \ref{relation inv - dexc}), 
$\hat{\pi}=\underline{4}\:1\:\underline{6}\:2\:\underline{8}\:\underline{9}
\:3\:5\:7\:10$ has as few inversions as possible for a permutation in ${\cal S}_{10}$ 
with excedance set $I$ and excedance letter set $A$, namely $12$.
\item We have both $\inv(\pi^{-1})=\inv(\pi)$ and $\dexc(\pi^{-1})=\dexc(\pi)$ for all $\pi\in{\cal S}_n$ where $\pi^{-1}$ 
denotes the inverse permutation. Furthermore, $\exc(\pi^{-1})=n-\exc(\pi)-{\sf 
fix}(\pi)$ for each $\pi\in{\cal S}_n$ where ${\sf fix}$ counts the number of fixed points in a permutation. Thus we 
can strengthen the upper bound in \ref{bounds} by replacing $\exc(\pi)$ with  
$\max\{\exc(\pi),n-\exc(\pi)-{\sf fix}(\pi)\}$.
\end{enum}
\erm
\end{rems}

For the investigations we have done until now, a certain kind of permutations 
plays an important role: those for which the inversion number and the excedance 
difference coincide. The following sections will deal with these, so-called 
bi-increasing, permutations in view of their excedances and descents. Here the 
idea is to construct one-to-one correspondences between bi-increasing permutations and 
other well-known combinatorial objects which encode the statistics as natural 
parameters.\\
In Section 2, first several simple equivalent definitions of bi-increasing permutations are given. It 
will turn out that these permutations are very interesting for everyone who 
studies forbidden patterns in permutations; they are just the $321$-avoiding 
ones. Moreover, their fixed points will be characterized.\\
Section 3 pursues the aim to determine the distribution of $\exc$ and $\dexc$ 
on the set of bi-increasing permutations. Recall that these permutations represent classes on 
which both statistics are constant. To this end, we establish first a 
one-to-one correspondence between bi-increasing permutations and step parallelogram polyominoes that transfers 
excedance number and excedance difference into width and area of the 
corresponding polyomino. By a simple transformation, the bijection can be given 
in terms of general parallelogram polyominoes as well. It proves that this 
correspondence connects the bijection between bi-increasing permutations and 
Dyck paths due to Billey, Jockusch, and Stanley with that one between 
parallelogram polyominoes and Dyck paths given by Delest and Viennot.\\
In Section 4, we develop an analogous model basing on 2-Motzkin paths to study 
the descent structure of bi-increasing permutations. First we will confine 
ourselves to such bi-increasing permutations whose excedances coincine with the 
descents. Here a bijection due to Fran\c{c}on and Viennot yields the desired 
correspondence where $\des$ and $\ddes$ are translated into the number of 
up-steps and the sum of height, respectively. By a little refinement we can 
extend this correspondence to all bi-increasing permutations. Moreover, the 
relation with another bijection between permutations and 2-Motzkin paths given 
by Foata and Zeilberger will be revealed.\\
Section 5 starts with the proof of the equidistribution of the difference 
statistics over the set of bi-increasing permutations. From this we obtain a 
one-to-one correspondence between parallelogram polyominoes and 2-Motzkin paths 
which transfers all the natural parameters to each other. As an application, we 
enumerate the skew diagrams (or, parallelogram polyominoes) according to their 
rank.\\
In the last section, we deduce a result concerning the 
symmetry of the joint distribution of $\exc$ and $\inv$ over ${\cal S}_n$ from 
an analogous result for bi-increasing permutations.
\vspace*{0.75cm}


\setcounter{section}{2}\setcounter{theo}{0}

\centerline{\large{\bf 2}\hspace*{0.25cm}
{\sc Characterization of bi-increasing permutations}}
\vspace*{0.5cm}

Theorem \ref{relation inv - dexc} characterizes the permutations for which inversion number 
and excedance difference are equally as those whose restrictions on excedances and non-excedances, 
respectively, are increasing words. In view of this, we call the elements of 
\bdpm
{\cal B}_n:=\{\pi\in{\cal S}_n:\inv(\pi)=\dexc(\pi)\}
\edpm
the {\it bi-increasing permutations} of length $n$. (This term was already used in \cite{foata-zeilberger}).
In particular, any bi-increasing permutation is uniquely determined by its excedances and 
excedance letters.\\[2ex]
In the following we give several simple equivalent definitions of ${\cal B}_n$.

\begin{prop} \label{bi-increasing via letters}
A permutation $\pi\in{\cal S}_n$ is bi-increasing if and only if 
\beas
&&\pi_k-k=|\{i>k:\pi_i<\pi_k\}|\quad\mbox{for }k\in{\sf E}(\pi),\\
&&k-\pi_k=|\{i<k:\pi_i>\pi_k\}|\quad\mbox{for }k\notin{\sf E}(\pi).
\eeas
\end{prop}

\begin{bew}
For any $k\in[n]$, the number of integers $i\in[n]$ with $\pi_i<\pi_k$ equals   
\bdpm
k-1-|\{i:i<k,\:\pi_i>\pi_k\}|+|\{i:i>k,\:\pi_i<\pi_k\}|.
\edpm
Hence $\pi_k=k+|\{i:i>k,\:\pi_i<\pi_k\}|-|\{i:i<k,\:\pi_i>\pi_k\}|$. If $k$ is 
an excedance then each integer $i<k$ with $\pi_i>\pi_k$ is excedance as well.
Analogously, every integer $i>k$ with $\pi_i<\pi_k$ is a non-excedance
if $k$ is such a one. In case of ordered words $\pi_{\sf e}$ and $\pi_{\sf ne}$, 
the sets $\{i:i<k,\:\pi_i>\pi_k\}$ for $k\in{\sf E}(\pi)$ and 
$\{i:i>k,\:\pi_i<\pi_k\}$ for $k\notin{\sf E}(\pi)$, respectively, are empty 
sets.
\end{bew}
\vspace*{1ex}

For $i=1,\ldots,n-1$, let $s_i$ denote the adjacent transposition $(i,i+1)\in{\cal 
S}_n$. It is well-known that every permutation $\pi\in{\cal S}_n$ 
can be written as $\pi=s_{i_1}s_{i_2}\cdots s_{i_k}$ where $k=\inv(\pi)$. The 
factorization of bi-increasing permutations is striking.

\begin{prop} \label{bi-increasing via factorization}
A permutation $\pi\in{\cal S}_n$ is bi-increasing if and only if
\bdpm
\pi=(s_{i_1+j_1-1}s_{i_1+j_1-2}\cdots s_{i_1})
(s_{i_2+j_2-1}s_{i_2+j_2-2}\cdots s_{i_2})\cdots 
(s_{i_e+j_e-1}s_{i_e+j_e-2}\cdots s_{i_e})
\edpm
for some positive integers $i_1,\ldots,i_e$ and $j_1,\ldots,j_e$ satisfying $1\le 
i_1<i_2<\ldots<i_e<n$ and $1<i_1+j_1<i_2+j_2<\ldots<i_e+j_e\le n$. In 
particular, $i_1,\ldots,i_e$ are precisely the excedances.
\end{prop}

\begin{bew}
We consider the decomposition of $\pi\in{\cal B}_n$ obtained by the following procedure:
\begin{enum1}
\item Set $k:=1$ and $\tau:=\pi$.
\item For $i=1,\ldots,n-1$, if $\tau_i>i$ then define $a_k:=\tau_i-1$ and 
$\tau:=s_{a_k}\tau$. Increase $k$ by $1$.
\item Write $\pi=s_{a_1}s_{a_2}\cdots s_{a_r}$ where $r$ denotes the value of 
$k$ before stopping (2).  
\end{enum1}
\vspace*{-0.3cm}

Let $i_1<\ldots<i_e$ be the excedances of $\pi$. The first step yields $a_1=\pi_{i_1}-1$. 
After the run for $i=i_1$, we obtain $a=(\pi_{i_1}-1,\pi_{i_1}-2,\ldots,i_1)$ and 
$\tau=s_{i_1}\cdots s_{\pi_{i_1}-2}s_{\pi_{i_1}-1}\pi$ where $1,2,\ldots,i_1$ are fixed 
points of $\tau$. Since $\pi_{\sf e}$ is increasing applying $s_{i_1}\cdots 
s_{\pi_{i_1}-1}$ to $\pi$ do not change the letters on the positions 
$i_2,\ldots,i_e$. Hence $\tau$ differs from $\pi$ in $\pi_{i_1}-i_1+1$ positions 
$j$, all different from $i_2,\ldots,i_e$. (We have $\tau_j=\pi_j+1$, except for 
$j=i_1$.) Note that $j$ is an non-excedance of $\tau$ as well. Assuming the contrary, $j>i_1$ has to be a 
fixed point of $\pi$. By Proposition \ref{bi-increasing via letters}, then there 
appears no letter greater than $j$ at a position to the left of $j$. Hence 
$j+1$ is not an excedance letter, and $j$ does not belong to the positions 
changed by the procedure. Consequently, ${\sf E}(\tau)={\sf 
E}(\pi)\setminus\{i_1\}$ and $\dexc(\tau)=\dexc(\pi)-(\pi_{i_1}-i_1)$. Applying 
the procedure successively to $\pi$ yields the desired factorization.\\
By reasoning similar, one shows that the inversion number of a permutation 
$\pi\in{\cal S}_n$ corresponding to the above product equals $\dexc(\pi)$.
\end{bew}

The following characterization deals with forbidden patterns and was already given in \cite[Lem. 5.6]{simion}. 
We say that a permutation $\pi\in{\cal S}_n$ {\it avoids the pattern} $321$ 
if there are no integers $i<j<k$ such that $\pi_i>\pi_j>\pi_k$, i.e., every 
decreasing subsequence in $\pi$ is of length at most two. (It is usual to write ${\cal S}_n(321)$ to denote 
the 321-avoiding permutations of length $n$.)

\begin{prop}
A permutation $\pi\in{\cal S}_n$ is bi-increasing if and only if it avoids the 
pattern $321$.
\end{prop}

\begin{bew}
The ``if" direction is evident. For the converse, first let $i<j$ be excedances. Obviously, 
there is an integer $k>j$ with $\pi_k\le j$. Thus we have $\pi_i<\pi_j$ since 
$\pi_i\pi_j\pi_k$ is decreasing otherwise. Clearly, the inverse $\pi^{-1}$ 
avoids $321$ if and only if $\pi$ does it. If $i<j$ are non-excedances but no 
fixed points then we obtain $\pi_i<\pi_j$ applying above argument to $\pi^{-1}$. Let 
now $i=\pi_i$ and $j>i$ be a non-excedance but not fixed. Assume that $i>\pi_j$. 
Clearly, then there exists an integer $k<i$ with $\pi_k>i$, that means, a
decreasing subsequence $\pi_k\,i\,\pi_j$. (In the reversed case, $i$ is 
non-excedance with $\pi_i<i$ and $j>i$ is fixed, we have $\pi_i<\pi_j$ anyway.)
\end{bew}

\begin{exam}
\brm
The permutation $2\:6\:1\:3\:7\:4\:5\:8\:10\:9\in{\cal S}_{10}$ will appear as 
example throughout the following sections. It is bi-increasing since 
$\inv(\pi)=\dexc(\pi)=8$, or equivalently, both $\pi_{\sf e}=2\:6\:7\:10$ and $\pi_{\sf 
ne}=1\:3\:4\:5\:8\:9$ are increasing words, or equivalently, 
$\pi=s_1s_5s_4s_3s_2s_6s_5s_9$, or 
equivalently, the maximum length of a decreasing subsequence of $\pi$ equals two. 
\erm
\end{exam}

Making use of these descriptions, we obtain a simple criterion for being a 
fixed point in a bi-increasing permutation.

\begin{cor} \label{fixed points}
An integer $i$ is a fixed point of $\pi\in{\cal B}_n$ 
if and only if $\pi_j<\pi_i$ for all $j<i$ and $\pi_j>\pi_i$ for all $j>i$.
\end{cor}

\begin{bew}
Let $i$ satisfy $\pi_i=i$ where $\pi\in{\cal B}_n$. By Proposition \ref{bi-increasing via 
letters}, there is no integer $j<i$ with $\pi_j>\pi_i$. An integer $j>i$ for which 
$\pi_j<\pi_i$ is a non-excedance which contradicts the fact that $\pi_{\sf ne}$ 
is increasing. The converse immediately follows from Proposition \ref{bi-increasing via 
letters}.
\end{bew}

In \cite[Th. 7.5]{robertson etal}, the authors enumerated bi-increasing 
permutations according to the number of their fixed points. There are 
\bdpm
\sum_{i=0}^{n-k} (-1)^i\left(\frac{k+1+i}{n+1}\right){2n-k-i\choose 
n}{k+i\choose k}
\edpm
permutations in ${\cal B}_n$ having exactly $k$ fixed points. In particular 
(see \cite[Cor. 3.3]{robertson etal}), the number of bi-increasing derangements of length $n$ is just the $n$th {\it Fine number} 
$F_n$ that may be defined by the formula $F_{n-1}+2F_n=C_n$. Here $C_n$ denotes 
the $n$th {\it Catalan number} defined by $\frac{1}{n+1}{2n\choose n}$. 

\begin{cor}
The number of bi-increasing permutations of length $n$ whose fixed points are 
exactly $i_1<i_2<\ldots<i_s$ equals $F_{i_1-1}F_{i_2-i_1-1}F_{i_3-i_2-1}\cdots 
F_{i_s-i_{s-1}-1}F_{n-i_s}$.
\end{cor}

\begin{bew}
By the previous corollary, any permutation 
$\pi\in{\cal B}_n$ having fixed points $i_1,\ldots,i_s$ can be represented as 
$\pi=\sigma_1i_1\sigma_2i_2\cdots\sigma_si_s\sigma_{s+1}$ where 
$\sigma_k$ is a bi-increasing derangement of length $i_k-i_{k-1}-1$ (set 
$i_0:=0$ and $i_{s+1}:=n+1$).
\end{bew}
\vspace*{0.75cm}


\setcounter{section}{3}\setcounter{theo}{0}

\centerline{\large{\bf 3}\hspace*{0.25cm}
{\sc Bi-increasing permutations and parallelogram polyominoes}}
\vspace*{0.5cm}

The characterization \ref{bi-increasing via factorization} inspires the 
following graphical representation of a bi-increasing permutation: for every 
excedance $i$ of $\pi\in{\cal B}_n$ draw a horizontal line of length $\pi_i-i$ beginning 
at $i$ such that the lines are arranged one below the other according to the 
appearance of $i$ in $\pi$.
\addtocounter{pic}{1}
\begin{center}
\unitlength0.4cm
\begin{picture}(9,5)
\linethickness{0.3pt}
\put(0,1){\line(1,0){9}}
\multiput(0,0.9)(1,0){10}{\line(0,1){0.2}}
\put(0,0.5){\makebox(0,0)[cc]{\tiny1}}
\put(1,0.5){\makebox(0,0)[cc]{\tiny2}}
\put(9,0.5){\makebox(0,0)[cc]{\tiny10}}
\bezier{30}(0,1)(0,3)(0,5)
\bezier{30}(9,1)(9,3)(9,5)
\put(0,5){\circle*{0.2}}\put(1,5){\circle*{0.2}}
\put(1,4){\circle*{0.2}}\put(5,4){\circle*{0.2}}
\put(4,3){\circle*{0.2}}\put(6,3){\circle*{0.2}}
\put(8,2){\circle*{0.2}}\put(9,2){\circle*{0.2}}
\linethickness{1pt}
\put(0,5){\line(1,0){1}}
\put(1,4){\line(1,0){4}}
\put(4,3){\line(1,0){2}}
\put(8,2){\line(1,0){1}}
\end{picture}

{\footnotesize{\bf Figure \thepic}\hspace*{0.25cm}Graphical representation of 
$2\:6\:1\:3\:7\:4\:5\:8\:10\:9\in{\cal B}_{10}$}
\end{center}

Therefore the distribution of the excedance number and excedance difference on 
${\cal B}_n$ answers the elementary question: in how many ways can a given 
number (corresponding to $\exc$) of lines of a prescribed total length (corresponding 
to $\dexc$) be arranged inside a stripe of width $n-1$ such that every line begins strictly to 
the right of the beginning of the previous line and ends strictly to the right 
of the end of the previous line?\\[2ex]
This representation immediately yields a connection between bi-increasing 
permutations and another class of well-known combinatorial objects, the 
parallelogram polyominoes. A {\it parallelogram polyomino} is a finite 
connected union of cells in the plane $\nat_0^2$ that can be described as 
region bounded by two non-intersecting lattice paths starting at the 
origin, using only steps $[0,1]$ and $[1,0]$, and ending in a common point. 
The most often studied parameters of parallelogram polyominoes 
are the {\it perimeter} which is the border length (twice path length), the 
{\it width} and {\it height} which are the coordinates of the point at that 
the paths end, and the {\it area} which is the number of cells.\\
For example, Figure 2 shows a parallelogram polyomino of perimeter 22, width 5, 
height 6, and area 13: 
\addtocounter{pic}{1}
\begin{center}
\unitlength0.4cm
\begin{picture}(5,5)
\linethickness{0.5pt}
\put(0,0){\line(1,0){1}}
\multiput(0,1)(0,1){2}{\line(1,0){2}}
\multiput(0,3)(0,1){2}{\line(1,0){4}}
\put(2,5){\line(1,0){3}}
\put(3,6){\line(1,0){2}}
\multiput(0,0)(1,0){2}{\line(0,1){4}}
\put(2,1){\line(0,1){4}}
\multiput(3,3)(1,0){2}{\line(0,1){3}}
\put(5,5){\line(0,1){1}}
\end{picture}

{\footnotesize{\bf Figure \thepic}\hspace*{0.25cm}A parallelogram polyomino}
\end{center}
                                 
First, we will consider a special kind of parallelogram polyominoes. We call a 
parallelogram polyomino a {\it step polyomino} if each horizontal border segment is of length 
$1$.

\begin{theo} \label{permutation-polyomino}
There is a bijection between bi-increasing permutations of length $n$ having 
$e$ excedances and excedance difference $k$ and step polyominoes of width $e+1$, height 
$n+1$, and area $k+n+1$. 
\end{theo}

\begin{bew}
Let $\pi\in{\cal B}_n$ be represented by a line arrangement as described above. 
Then the $i$th line (counted from left) is defined to correspond to the rows in common 
to the columns $i$ and $i+1$ of the step polyomino. Let the first column begin at level 
$0$, and let the last column end at level $n+1$. It is easy to see that the 
cell number of the resulting polyomino is just the total length of all lines plus $n+1$.
\end{bew}

\begin{exam}
\brm
Let $\pi=2\:6\:1\:3\:7\:4\:5\:8\:10\:9\in{\cal B}_{10}$ as before. The 
corresponding step polyomino is:
\addtocounter{pic}{1}
\begin{center}
\unitlength0.37cm
\begin{picture}(6,11)
\linethickness{0.3pt}
\put(1,0){\line(0,1){11}}
\multiput(0.9,0)(0,1){12}{\line(1,0){0.2}}
\put(0.5,0){\makebox(0,0)[cc]{\tiny0}}
\put(0.5,1){\makebox(0,0)[cc]{\tiny1}}
\put(0.5,11){\makebox(0,0)[cc]{\tiny11}}
\linethickness{0.5pt}
\put(1,0){\line(1,0){1}}
\put(1,1){\line(1,0){2}}
\put(1,2){\line(1,0){3}}
\multiput(2,3)(0,1){2}{\line(1,0){2}}
\multiput(2,5)(0,1){2}{\line(1,0){3}}
\put(3,7){\line(1,0){2}}
\put(4,8){\line(1,0){1}}
\multiput(4,9)(0,1){2}{\line(1,0){2}}
\put(5,11){\line(1,0){1}}
\put(1,0){\line(0,1){2}}
\multiput(2,0)(1,1){2}{\line(0,1){6}}
\put(4,2){\line(0,1){8}}
\put(5,5){\line(0,1){6}}
\put(6,9){\line(0,1){2}}
\linethickness{1pt}
\put(2,1){\line(0,1){1}}
\put(3,2){\line(0,1){4}}
\put(4,5){\line(0,1){2}}
\put(5,9){\line(0,1){1}}
\end{picture}

{\footnotesize{\bf Figure \thepic}\hspace*{0.25cm}Step polyomino associated 
with $2\:6\:1\:3\:7\:4\:5\:8\:10\:9$}
\end{center}
\erm
\end{exam}

For a more formal description, we code a step polyomino by two integer 
sequences that contain the lengths of the vertical border segments. Let 
$\alpha=(l_1-1,l_2,l_3,\ldots,l_w)$ and $\beta=(k_1,k_2,\ldots,k_{w-1},k_w-1)$ 
where $l_i$ and $k_i$ denote the lengths of the $i$th vertical border segment of the 
upper and lower path, respectively. We will identify a step polyomino with the 
pair $(\alpha,\beta)$ and vice versa.\\
For instance, Figure 3 shows the step polyomino $((1,4,1,3,1),(1,1,3,4,1))$.\\
Given a step polyomino of width $w$ and height $n+1$, obviously, $\alpha$ and $\beta$ are compositions of $n$ 
into $w$ positive parts. On the other hand, it is easy to say when a 
pair of compositions describes a step polyomino.

\begin{prop} \label{dominance condition}
Let ${\cal C}_{n,w}$ be the set of compositions of $n$ into $w$ positive parts.
Any pair $(\alpha,\beta)\in{\cal C}_{n,w}^2$ corresponds to a step polyomino of width 
$w$ and height $n+1$ if and only if $\alpha\ge\beta$ (dominance order), i.e.,   
$\alpha_1+\ldots+\alpha_i\ge\beta_1+\ldots+\beta_i$ for all $i\ge1$.   
\end{prop}

\begin{bew}
If the dominance condition fails for some integer $i$ then the border 
paths intersect at the point $(i,\alpha_1+\ldots+\alpha_i+1)$.
\end{bew}
\vspace*{1ex}

To determine the distribution of $\exc$ on bi-increasing permutations first we 
give a bijection between step polyominoes and objects whose numbers is well-known.

\begin{theo} \label{polyomino-diagram}
There is a bijection between step polyominoes of width $w$ and height $n+1$ 
and Young diagrams which fit into the shape $(n-1,n-2,\ldots,1)$ and have 
$w-1$ corners. 
\end{theo}

\begin{bew}
Let $(\alpha,\beta)\in{\cal C}_{n,w}^2$. Define 
$\lambda=(n-\alpha_1,n-\alpha_1-\alpha_2,\ldots,n-\alpha_1-\ldots-\alpha_{w-1})$ 
and $\mu=(\beta_1+\ldots+\beta_{w-1},\ldots,\beta_1+\beta_2,\beta_1)$. Hence 
$\lambda$ and $\mu$ are partitions into distinct parts, all at most $n-1$. (The Young 
diagrams of $\lambda$ and $\mu$ appear when we consider the polyomino as being 
contained in an $w\times(n+1)$-rectangle.)\\
There is exactly one partition $\Lambda$ such that 
\begin{enum2}
\item the {\it distinct} parts of $\Lambda$ are just the parts of $\lambda$, and
\item the {\it distinct} parts of the conjugate to $\Lambda$ are just the parts of 
$\mu$.
\end{enum2}
\vspace*{-0.35cm}

Given $\lambda$ and $\mu$, we construct $\Lambda$ as follows. Beginning with 
the diagram of $\lambda$, add $\mu_1-(w-1)$ squares to each of the first $\lambda_{w-1}$ columns, then add $\mu_2-(w-2)$ squares to each of the next $\lambda_{w-2}-\lambda_{w-1}$ 
columns, then add $\mu_3-(w-3)$ squares to each of the next 
$\lambda_{w-3}-\lambda_{w-2}$ columns, and so on. Since $\mu_{w-i}\ge i$ 
this procedure is always possible. By the construction, the 
squares $(\mu_{w-i},\lambda_i)$ where $i=1,\ldots,w-1$ are just the corners of 
the Young diagram of $\Lambda$. (Here the first coordinate denotes the rows, and 
the second one denotes the columns of the diagram, beginning at the left-hand 
upper corner.) Thus the partition 
$\Lambda$ is uniquely determined, and it satisfies (i) and (ii).\\ 
In addition, Proposition \ref{dominance condition} yields $\mu_{w-i}+\lambda_i=\beta_1+\ldots+\beta_i+n-\alpha_1-\ldots-\alpha_i\le 
n$ for each $i$, that is, $\Lambda$'s diagram fits into $(n-1,n-2,\ldots,1)$. 
\end{bew}

\begin{exam}
\brm
For the step polyomino $((1,4,1,3,1),(1,1,3,4,1))$ we obtain the partitions 
$\lambda=(9,5,4,1)$ and $\mu=(9,5,2,1)$, and hence $\Lambda=9\:5\:4^3\:1^4$.
\addtocounter{pic}{1}
\begin{center}                                                        
\definecolor{gray1}{gray}{0.85}
\fboxsep0cm
\fboxrule0cm
\unitlength0.4cm
\begin{picture}(22,12)
\put(0,2){\fcolorbox{gray1}{gray1}{\makebox(1,9){}}}
\put(1,6){\fcolorbox{gray1}{gray1}{\makebox(1,5){}}}
\put(2,7){\fcolorbox{gray1}{gray1}{\makebox(1,4){}}}
\put(3,10){\fcolorbox{gray1}{gray1}{\makebox(1,1){}}}
\put(13,10){\fcolorbox{gray1}{gray1}{\makebox(9,1){}}}
\put(13,9){\fcolorbox{gray1}{gray1}{\makebox(5,1){}}}
\put(13,8){\fcolorbox{gray1}{gray1}{\makebox(4,1){}}}
\put(13,7){\fcolorbox{gray1}{gray1}{\makebox(1,1){}}}
\linethickness{0.5pt}
\put(0,0){\line(1,0){5}}
\put(0,1){\line(1,0){2}}
\put(0,2){\line(1,0){3}}
\multiput(1,3)(0,1){2}{\line(1,0){2}}
\multiput(1,5)(0,1){2}{\line(1,0){3}}
\put(2,7){\line(1,0){2}}
\put(3,8){\line(1,0){1}}
\multiput(3,9)(0,1){2}{\line(1,0){2}}
\put(0,11){\line(1,0){5}}
\put(0,0){\line(0,1){11}}
\multiput(1,0)(1,1){2}{\line(0,1){6}}
\put(3,2){\line(0,1){8}}
\put(4,5){\line(0,1){6}}
\put(5,0){\line(0,1){11}}
\multiput(13,10)(0,1){2}{\line(1,0){9}}
\put(13,9){\line(1,0){5}}
\multiput(13,6)(0,1){3}{\line(1,0){4}}
\multiput(13,2)(0,1){4}{\line(1,0){1}}
\multiput(13,2)(1,0){2}{\line(0,1){9}}
\multiput(15,6)(1,0){3}{\line(0,1){5}}
\put(18,9){\line(0,1){2}}
\multiput(19,10)(1,0){4}{\line(0,1){1}}
\put(9,6){\makebox(0,0)[cc]{$\longleftrightarrow$}}
\end{picture}

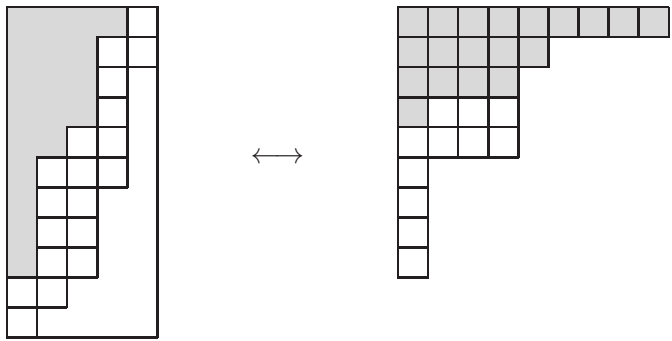
{\footnotesize{\bf Figure \thepic}\hspace*{0.25cm}One-to-one correspondence between step 
polyomino and restricted Young diagram}
\end{center}
\erm
\end{exam}

\begin{rem} \label{partial sums}
\brm
Connecting the bijections given in \ref{permutation-polyomino} and 
\ref{polyomino-diagram} yields a one-to-one correspondence between bi-increasing 
permutations and restricted Young diagrams.\\ 
Let $(\alpha,\beta)\in{\cal C}_{n,w}^2$ be the step polyomino associated with $\pi\in{\cal B}_n$. 
Then the partial sums $\beta_1,\beta_1+\beta_2,\ldots,\beta_1+\ldots+\beta_{w-1}$ 
and $\alpha_1+1,\alpha_1+\alpha_2+1,\ldots,\alpha_1+\ldots+\alpha_{w-1}+1$ are exactly the excedances and excedance letters of $\pi$, respectively. 
Therefore the bi-increasing permutation $\pi$ corresponds to the Young diagram 
having corners $(i_k,n+1-\pi_{i_k})$ where $i_k$ denotes the excedances of 
$\pi$. 
\erm
\end{rem}

The restricted Young diagrams are an item in Stanley's list of combinatorial objects counted by 
Catalan numbers (see \cite[Ex. 6.19]{stanley1}).\\ 
Consider the diagram as being contained in an $n\times n$-rectangle. Then the 
lattice path going along the diagram boundary from the lower left-hand 
rectangle corner to the upper right-hand one never falls below the line $x=y$. Paths of this kind 
were counted in \cite{sulanke}. There are 
\bdpm
N_{n,w}=\frac{1}{n}{n\choose w}{n\choose w-1}
\edpm
paths with $w$ horizontal segments, or equivalently, Young diagrams 
with $w-1$ corners fitting into $(n-1,n-2,\ldots,1)$. The integers $N_{n,w}$ 
are called {\it Narayana numbers}.

\begin{cor} \label{exc-distribution}
The statistics $\exc$ is Narayana distributed over ${\cal B}_n$; we have 
\bdpm
B_n^{\exc}(e)=\frac{1}{n}{n\choose e}{n\choose e+1}\quad\mbox{for all $n$ and 
$e$}.
\edpm
\end{cor}

Making use of the several bijections established above, we obtain 
the following combinatorial interpretations of Catalan numbers and Narayana 
numbers, respectively. 

\begin{cor} \label{Catalan interpretations}
The Catalan number $C_n$ (Narayana number $N_{n,w}$) counts the number of 
\begin{enum}
\item arrangements of ($w-1$) lines inside a stripe of width $n-1$, as described above,
\item step polyominoes of height $n+1$ (and width $w$),
\item pairs $(\alpha,\beta)$ of compositions of $n$ with the same number ($w$) of 
parts such that $\alpha\ge\beta$. 
\end{enum}
\end{cor}
\vspace*{-5ex}

\begin{rem}
\brm
Certainly, the bijection given in \ref{polyomino-diagram} is not the most obvious one proving 
Corollary \ref{exc-distribution}. In the following we give an elementary transformation that  
takes step polyominoes to general parallelogram polyominoes, and proves \ref{exc-distribution} as 
well. However, the one-to-one correspondence between bi-increasing (or $321$-avoiding) 
permutations and restricted Young diagrams yields a simple bijection between 
those permutations and $132$-avoiding ones which we have described in 
\cite{reifegerste}.      
\erm
\end{rem}

Let $(\alpha,\beta)\in{\cal C}_{n,w}^2$ be a step polyomino of width $w$ and 
height $n+1$. Define the components of the sequences 
\bdpm
\gamma=(\alpha_1,\alpha_2-1,\alpha_3-1,\ldots,\alpha_w-1)\quad\mbox{and}\quad
\delta=(\beta_1-1,\beta_2-1,\ldots,\beta_{w-1}-1,\beta_w)
\edpm
to be the vertical segment lengths of the lattice paths bounding a parallelogram 
polyomino having perimeter $\sum_{i=1}^w (\gamma_i+\delta_i)+2w=2n+2$. (The 
border paths are generated by adjoining alternately $\gamma_i$ and $\delta_i$, 
respectively, north steps and one east step.) Since $\alpha\ge\beta$ the path 
described by $\delta$ never rises above the path encoded by $\gamma$. Note that $\gamma$ and $\delta$ are compositions of 
$n+1-w$ into {\it non-negative} parts. Hence the resulting polyomino is of 
width $w$ and height $n+1-w$.\\ 
For example, the step polyomino $((1,4,1,3,1),(1,1,3,4,1))$ is transformed into 
the parallelogram polyomino $((1,3,0,2,0),(0,0,2,3,1))$:
\addtocounter{pic}{1}
\begin{center}                                                        
\unitlength0.4cm
\begin{picture}(22,12)
\linethickness{0.5pt}
\put(0,0){\line(1,0){1}}
\put(0,1){\line(1,0){2}}
\put(0,2){\line(1,0){3}}
\multiput(1,3)(0,1){2}{\line(1,0){2}}
\multiput(1,5)(0,1){2}{\line(1,0){3}}
\put(2,7){\line(1,0){2}}
\put(3,8){\line(1,0){1}}
\multiput(3,9)(0,1){2}{\line(1,0){2}}
\put(4,11){\line(1,0){1}}
\put(0,0){\line(0,1){2}}
\multiput(1,0)(1,1){2}{\line(0,1){6}}
\put(3,2){\line(0,1){8}}
\put(4,5){\line(0,1){6}}
\put(5,9){\line(0,1){2}}
\multiput(13,0)(0,1){2}{\line(1,0){3}}
\multiput(14,2)(0,1){3}{\line(1,0){3}}
\multiput(16,5)(0,1){2}{\line(1,0){2}}
\put(13,0){\line(0,1){1}}
\multiput(14,0)(1,0){2}{\line(0,1){4}}
\put(16,0){\line(0,1){6}}
\put(17,2){\line(0,1){4}}
\put(18,5){\line(0,1){1}}
\put(9,4){\makebox(0,0)[cc]{$\longrightarrow$}}
\end{picture}

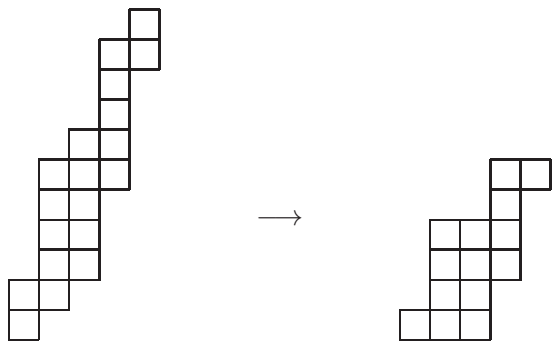
{\footnotesize{\bf Figure \thepic}\hspace*{0.25cm}Transformation of a step polyomino into 
a general parallelogram polyomino}
\end{center}
It is known (for instance, see \cite{delest}) that the number of parallelogram polyominoes 
having width $w$ and perimeter $2n+2$ equals the Narayana number $N_{n,w}$.\\
Obviously, the area of $(\gamma,\delta)$ is equal to the area of 
$(\alpha,\beta)$ minus $w$. Thus parallelogram polyominoes encode $\exc$ and $\dexc$ 
by their width and area. More exactly we have

\begin{theo} \label{permutation-parallelogram polyomino}
Let $\pi\in{\cal B}_n$ be a bi-increasing permutation with $\exc(\pi)=e$ and 
$\dexc(\pi)=k$. Then $\pi$ corresponds in a one-to-one fashion to a 
parallelogram polyomino of perimeter $2n+2$, width $e+1$, and area $n+k-e$.
\end{theo}

In \cite{barcucci etal2}, the authors gave a bijection between $321$-avoiding 
permutations of length $n$ and parallelogram polyominoes of perimeter $2n+2$. 
Their idea is as follows.\\
West described in \cite{west} a method to construct the set ${\cal B}_n$ 
recursively, by means of generating trees. A basic term in this context is that 
of active sites of a permutation. Given $\pi\in{\cal B}_n$, an integer 
$i\in[n+1]$ is called an {\it active site} of $\pi$ if the permutation 
$\pi^{(i)}:=\pi_1\pi_2\cdots\pi_{i-1}\:n+1\:\pi_i\pi_{i+1}\cdots\pi_n$ avoids $321$ 
as well. Barcucci et al. studied the consequences for the number of active sites and the inversion 
number when inserting $n+1$ into the site $i$.\\ 
On the side of parallelogram polyominoes, all polyominoes of perimeter $2(n+1)$ 
can be constructed from those having perimeter $2n$ by adding a cell onto the 
last column, or adding a column to the right of the last column such that both 
end at the same level. Here \cite{barcucci etal2} considered the changes of 
width, height, cell number of the last column, and number of cells 
having an adjacent cell on its right under these operations.\\
The comparison yields a correspondence between the following parameters:   
\beas
\mbox{\bf bi-increasing permutations}&&\mbox{\bf parallelogram polyominos}\\
\mbox{length }n&&\mbox{perimeter }2n+2\\
\mbox{number of active sites minus one}&&\mbox{number of cells belonging to the last 
column}\\
\mbox{inversion number}&&\mbox{number of cells with adjacent cell on its right}
\eeas

It is easy to see that our bijection carries out these transformations.

\begin{prop}
The bijection given in {\rm \ref{permutation-parallelogram polyomino}} 
translates the permutation statistics into the polyomino parameters, as 
described above.
\end{prop}

\begin{bew}
Let $\pi\in{\cal B}_n$ and $(\gamma,\delta)$ its corresponding parallelogram 
polyomino. It is evident from the structure of bi-increasing permutations that 
$\pi^{(i)}\in{\cal B}_{n+1}$ if (and only if) $i>i_e$ where $i_e$ denotes the 
greatest excedance of $\pi$. By construction, the last column of 
$(\gamma,\delta)$ contains precisely $n-i_e$ cells. Clearly, the number of 
cells having an adjacent cell on its right equals area minus height of 
$(\gamma,\delta)$. (Note that the half-perimeter is just the sum of height and 
width.) Since $\inv(\pi)=\dexc(\pi)$, the bijection \ref{permutation-parallelogram 
polyomino} translates the inversion number as desired.\\
In detail, $\pi^{(n+1)}$ corresponds to the parallelogram polyomino obtained from 
$(\gamma,\delta)$ by adding an addition cell on top of the last column. (The 
permutations $\pi$ and $\pi^{(n+1)}$ have the same excedances and excedance 
letters.) For $i=i_e+1,\ldots,n$, a new excedance with difference $n+1-i$ is 
arisen in $\pi^{(i)}$. (All the other excedances are preserved with their letters.) 
Thus $\pi^{(i)}$'s parallelogram polyomino develops from $(\gamma,\delta)$ by 
adding an additional column to the right of the last column, 
consisting of $n+1-i$ cells, each with an adjacent cell on its left. 
\end{bew}

\begin{exam}
\brm
Let $\pi=4\:1\:2\:5\:3\:6\in{\cal B}_6$, with maximal excedance $4$. Thus 
$\pi^{(5)}$, $\pi^{(6)}$, and $\pi^{(7)}$ are bi-increasing permutations 
corresponding to the following parallelogram polyominoes:
\addtocounter{pic}{1}
\begin{center}                                                        
\definecolor{gray1}{gray}{0.85}
\fboxsep0cm
\fboxrule0cm
\unitlength0.4cm
\begin{picture}(32,7)
\put(12,3.5){\fcolorbox{gray1}{gray1}{\makebox(1,2){}}}
\put(22,4.5){\fcolorbox{gray1}{gray1}{\makebox(1,1){}}}
\put(31,5.5){\fcolorbox{gray1}{gray1}{\makebox(1,1){}}}
\linethickness{0.5pt}
\multiput(0,1.5)(0,1){2}{\line(1,0){2}}
\multiput(0,3.5)(0,1){2}{\line(1,0){3}}
\put(2,5.5){\line(1,0){1}}
\multiput(0,1.5)(1,0){2}{\line(0,1){3}}
\put(2,1.5){\line(0,1){4}}
\put(3,3.5){\line(0,1){2}}
\multiput(9,1.5)(0,1){2}{\line(1,0){2}}
\multiput(9,3.5)(0,1){2}{\line(1,0){4}}
\put(11,5.5){\line(1,0){2}}
\multiput(9,1.5)(1,0){2}{\line(0,1){3}}
\put(11,1.5){\line(0,1){4}}
\multiput(12,3.5)(1,0){2}{\line(0,1){2}}
\multiput(19,1.5)(0,1){2}{\line(1,0){2}}
\put(19,3.5){\line(1,0){3}}
\put(19,4.5){\line(1,0){4}}
\put(21,5.5){\line(1,0){2}}
\multiput(19,1.5)(1,0){2}{\line(0,1){3}}
\put(21,1.5){\line(0,1){4}}
\put(22,3.5){\line(0,1){2}}
\put(23,4.5){\line(0,1){1}}
\multiput(29,1.5)(0,1){2}{\line(1,0){2}}
\multiput(29,3.5)(0,1){2}{\line(1,0){3}}
\multiput(31,5.5)(0,1){2}{\line(1,0){1}}
\multiput(29,1.5)(1,0){2}{\line(0,1){3}}
\put(31,1.5){\line(0,1){5}}
\put(32,3.5){\line(0,1){3}}
\put(1.5,0){\makebox(0,0)[cb]{\small$\pi=4\:1\:2\:5\:3\:6$}}
\put(11,0){\makebox(0,0)[cb]{\small$\pi^{(5)}=4\:1\:2\:5\:7\:3\:6$}}
\put(21,0){\makebox(0,0)[cb]{\small$\pi^{(6)}=4\:1\:2\:5\:3\:7\:6$}}
\put(30.5,0){\makebox(0,0)[cb]{\small$\pi^{(7)}=4\:1\:2\:5\:3\:6\:7$}}
\end{picture}

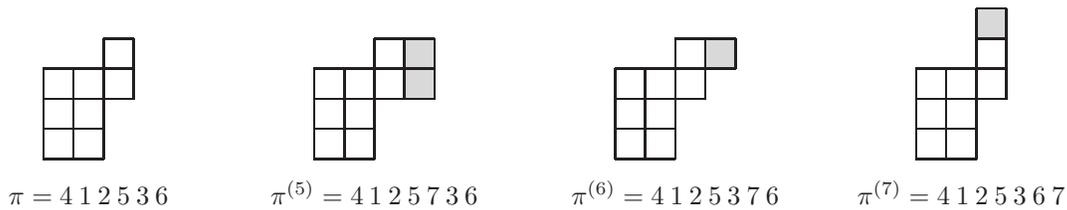
\parbox[t]{14cm}{{\footnotesize{\bf Figure \thepic}\hspace*{0.25cm}Permutations obtained from 
$4\:1\:2\:5\:3\:6$ by West's method with their corresponding\\[-1.5ex]
parallelogram polyominoes}}
\end{center}
\erm
\end{exam}
\vspace*{1ex}

Using a result of \cite{west}, the correspondence between the maximum excedance and the number of active sites 
in bi-increasing permutations yields the following enumerative statement.

\begin{cor}
There are ${n-1+k\choose k}-{n-1+k\choose k-1}$ bi-increasing permutations of length $n$ whose greatest 
excedance equals $k$.
\end{cor}

\begin{bew}
Clearly, the number of active sites of $\pi\in{\cal B}_n$ having greatest 
excedance $k$ is equal to $n+1-k$. By \cite[Th. 2.12]{west} (applied on reverse 
permutations), we obtain the given number.
\end{bew}

In \cite{barcucci etal2}, the authors also determined the generating function 
of ${\cal B}_n$ according to length and inversion number. This yields the 
distribution coefficient $B_n^{\dexc}(k)$ we have looked for. (As shown above, $B_n^{\dexc}(k)$ is the number of 
parallelogram polyominoes of perimeter $2n+2$ for which area minus height equals $k$.)

\begin{theo}[{\cite[Th. 3.4]{barcucci etal2}}] \label{series}
The number of bi-increasing permutations having length $n$ and excedance 
difference $k$ is the coefficient of $x^nq^k$ in the quotient 
\bdpm
\frac{J_1(x,q)}{J_0(x,q)}
\edpm
where the function $J_r$ is defined by
\bdpm
J_r(x,q)=\sum_{n\ge0} 
\frac{(-1)^nx^{n+r}q^{\frac{1}{2}n(n+2r+1)}}{(x)_{n+r}(q)_n}
\edpm 
with $(a)_n:=(1-a)(1-aq)(1-aq^2)\cdots(1-aq^{n-1})$.
\end{theo}

Because the enumeration of polyominoes is a topic on its own we leave it here 
at the combinatorial interpretation for the distribution coefficients 
$B_n^{(\exc,\dexc)}(e,k)$. We will determine these integers in a forthcoming 
paper.\\
A nice property of the joint distribution of $\exc$ and $\dexc$ on the set of 
bi-increasing permutations can be proved already now.

\begin{cor} \label{symmetry}
$B_n^{(\exc,\dexc)}(e,k)=B_n^{(\exc,\dexc)}(n-1-e,n-1-2e+k)$ for all $e$ and $k$.
\end{cor}

\begin{bew}
We give an elementary involution on the parallelogram polyominoes which yields 
the desired symmetry. Given a parallelogram polyomino, we reverse each of the border paths 
and replace a north step {\sf N} by an east step {\sf E} and vice versa.\\
For example, $((1,3,0,2,0),(0,0,2,3,1))$ corresponds to 
$((2,0,2,0,0,1),(1,0,0,1,0,3))$ in this way: 
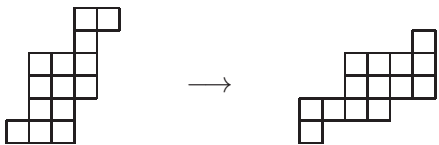
\begin{center}                                                        
\unitlength0.3cm
\begin{picture}(19,6)
\linethickness{0.5pt}
\multiput(0,0)(0,1){2}{\line(1,0){3}}
\multiput(1,2)(0,1){3}{\line(1,0){3}}
\multiput(3,5)(0,1){2}{\line(1,0){2}}
\put(0,0){\line(0,1){1}}
\multiput(1,0)(1,0){2}{\line(0,1){4}}
\put(3,0){\line(0,1){6}}
\put(4,2){\line(0,1){4}}
\put(5,5){\line(0,1){1}}
\put(13,0){\line(1,0){1}}
\put(13,1){\line(1,0){4}}
\put(13,2){\line(1,0){6}}
\multiput(15,3)(0,1){2}{\line(1,0){4}}
\put(18,5){\line(1,0){1}}
\multiput(13,0)(1,0){2}{\line(0,1){2}}
\multiput(15,1)(1,0){3}{\line(0,1){3}}
\multiput(18,2)(1,0){2}{\line(0,1){3}}
\put(9,2.5){\makebox(0,0)[cc]{$\longrightarrow$}}
\end{picture}
\end{center}
(The upper border path {\sf NENNNEENNEE} of the left-hand polyomino is transformed into {\sf 
NNEENNEEENE}, the upper border path of the right-hand polyomino.)\\ 
Obviously, perimeter and area are not changed by this map. The width (resp. 
height) of a polyomino equals the height (resp. width) of the corresponding 
one. Theorem \ref{permutation-parallelogram polyomino} yields the assertion.
\end{bew}
   
$321$-avoiding permutations have been studied in a manifold way. In particular, 
several authors gave one-to-one correspondences to lattice paths (for 
example, see \cite{billey etal} and \cite{krattenthaler}). Before we do this in 
the next section as well, we complete this section with a note on the close connection between two 
well-known bijections dealing with Dyck paths. A {\it Dyck path} is a lattice path in the plane $\nat_0^2$ 
from the origin to $(2n,0)$ consisting only of up-steps $[1,1]$ and down-steps 
$[1,-1]$. Any point of the Dyck path connecting an up-step with a following down-step we call a {\it 
peak}, any point connecting a down-step with a following up-step we call a {\it 
valley} of the path.\\
In \cite[p. 361]{billey etal}, Billey, Jockusch, and Stanley established the  
bijection $\Phi_{BJS}$ between $321$-avoiding permutations of length $n$ and 
Dyck paths of length $2n$.\\ 
In \cite[Sect. 4]{delest-viennot}, Delest and Viennot gave the bijection $\Phi_{DV}$ 
between parallelogram polyominoes of perimeter $2n+2$ and Dyck paths having 
length $2n$.\\
It is not difficult to see that the bijection from \ref{permutation-parallelogram polyomino} 
connects the both ones. 

\begin{prop}
Let $\pi\in{\cal B}_n$ be a bi-increasing permutation and $(\gamma,\delta)$ 
its corresponding parallelogram polyomino. Then we have 
$\Phi_{BJS}(\pi)=\Phi_{DV}((\gamma,\delta))$. 
\end{prop}

\begin{bew}
Let $i_1<\ldots<i_e$ be the excedances of $\pi\in{\cal B}_n$. The bijection $\Phi_{BJS}$ 
takes $\pi$ to the Dyck path defined as follows. For $k=1,\ldots,e$, let 
$a_k:=\pi_{i_k}-1$, and $a_0:=0,\;a_{e+1}:=n$. Furthermore, for $k=1,\ldots,e$, let 
$b_k:=i_k$, and $b_0:=0,\;b_{e+1}:=n$. Beginning at the origin adjoin 
alternately $a_k-a_{k-1}$ up-steps and $b_k-b_{k-1}$ down-steps where 
$k=1,\ldots,e+1$. Note that $a_k-a_{k-1}=\alpha_k$ and $b_k-b_{k-1}=\beta_k$ 
for $k=1,\ldots,e+1$. Here $(\alpha,\beta)$ is the step polyomino associated 
with the permutation $\pi$ by Theorem \ref{permutation-polyomino}.\\
By the construction, the heights (the ordinates) of the peaks are $a_k-b_{k-1}$ 
with $k=1,\ldots,e+1$. The valleys are of height $a_k-b_k$ where $k=1,\ldots,e$.\\
On the other hand, the number of cells belonging to the $k$th column of $(\gamma,\delta)$ is just 
\bdpm
\gamma_1+\ldots+\gamma_k-(\delta_1+\ldots+\delta_{k-1})=\alpha_1+\ldots+\alpha_k-
(\beta_1+\ldots+\beta_{k-1})=a_k-b_{k-1}.
\edpm
For the number of cells adjacent to the columns $k$ and $k+1$ of 
$(\gamma,\delta)$ we obtain
\bdpm
\gamma_1+\ldots+\gamma_k-(\delta_1+\ldots+\delta_k)=\alpha_1+\ldots+\alpha_k+1-
(\beta_1+\ldots+\beta_k)=a_k-b_k+1.
\edpm
The correspondence between column sizes and peak heights and between the numbers 
of rows in common of adjacent columns and valley heights, respectively, is 
precisely the description of $\Phi_{DV}$.
\end{bew}

By the way, the Young diagram constructed in the proof of Theorem 
\ref{polyomino-diagram} we also find again: the Dyck path describes just the boundary 
of the diagram when considered as being contained in a $n\times n$-rectangle.

\begin{exam}
\brm
We consider $\pi=2\:6\:1\:3\:7\:4\:5\:8\:10\:9\in{\cal B}_{10}$ again. Billey-Jockusch-Stanley's bijection takes $\pi$ to the Dyck path
\addtocounter{pic}{1}
\begin{center}                                                        
\unitlength0.4cm
\begin{picture}(20,4)
\linethickness{0.5pt}
\bezier{50}(0,0)(0.5,0.5)(1,1)
\bezier{50}(1,1)(1.5,0.5)(2,0)
\bezier{200}(2,0)(4,2)(6,4)
\bezier{50}(6,4)(6.5,3.5)(7,3)
\bezier{50}(7,3)(7.5,3.5)(8,4)
\bezier{150}(8,4)(9.5,2.5)(11,1)
\bezier{150}(11,1)(12.5,2.5)(14,4)
\bezier{200}(14,4)(16,2)(18,0)
\bezier{50}(18,0)(18.5,0.5)(19,1)
\bezier{50}(19,1)(19.5,0.5)(20,0)
\linethickness{0.2pt}
\multiput(0,0)(0,1){5}{\line(1,0){20}}
\multiput(0,0)(1,0){21}{\line(0,1){4}}
\put(1,1){\circle*{0.25}}
\put(6,4){\circle*{0.25}}
\put(8,4){\circle*{0.25}}
\put(14,4){\circle*{0.25}}
\put(19,1){\circle*{0.25}}
\put(2,0){\circle{0.25}}
\put(7,3){\circle{0.25}}
\put(11,1){\circle{0.25}}
\put(18,0){\circle{0.25}}
\end{picture}

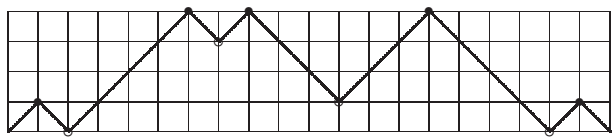
{\footnotesize{\bf Figure \thepic}\hspace*{0.25cm}Dyck path 
$\Phi_{BJS}(2\:6\:1\:3\:7\:4\:5\:8\:10\:9)$ resp. $\Phi_{DV}((1,3,0,2,0),(0,0,2,3,1))$}
\end{center}
which is exactly the path obtained from the parallelogram polyomino corresponding 
to $\pi$ (see the right-hand polyomino in Figure 5) by Delest-Viennot's bijection. Note that 
$\exc(\pi)$ is the number of valleys and $\dexc(\pi)$ equals the sum of 
the heights of the valleys, each increased by 1.
\erm
\end{exam}

\begin{rem} \label{derangement}
\brm
Corollary \ref{fixed points} says, if an integer $i$ is a fixed point of $\pi\in{\cal B}_n$ then there is no line in $\pi$'s graphical representation which 
starts or ends at level $i$, or runs through the point $i$. Consequently, the corresponding step polyomino has two consecutive rows consisting of one cell. 
The transformation into general parallelogram polyominoes deletes one of these cells. 
Therefore, the number of fixed points of a bi-increasing permutation equals the 
number of singleton rows of the associated parallelogram polyomino.\\
Hence there are $F_n$ parallelogram polyominoes of perimeter $2n+2$ whose all rows are of length at most 
two (see also \cite{deutsch-shapiro}'s list of objects counted by Fine 
numbers).
\erm
\end{rem}
\vspace*{0.6cm}


\setcounter{section}{4}\setcounter{theo}{0}

\centerline{\large{\bf 4}\hspace*{0.25cm}
{\sc Bi-increasing permutations and 2-Motzkin paths}}
\vspace*{0.5cm}

From the nature of bi-increasing permutations follows: an integer $i$ is a 
descent if and only if $i$ is an excedance but $i+1$ is none. In particular, 
there are no consecutive descents.\\
Consequently, the excedance number and descent number have not the same 
distribution over ${\cal B}_n$. We have $\des(\pi)\le\frac{n}{2}$ for all 
$\pi\in{\cal B}_n$.\\[2ex]
In this section we will develop a correspondence between bi-increasing 
permutations and certain lattice paths that makes it possible to determine the 
distribution of the descent number and descent difference on ${\cal B}_n$ by 
means of path enumeration.\\
First we consider only a subset of ${\cal B}_n$, namely those permutations with identical excedance number and 
descent number. Note that the condition $\exc(\pi)=\des(\pi)$ for $\pi\in{\cal 
B}_n$ is equivalent to ${\sf E}(\pi)={\sf D}(\pi)$. Permutations satisfying 
this can easily be characterized.

\begin{prop} \label{increasing descent and non-descent words}
Let $\pi\in{\cal S}_n$ be a permutation. Then the words $\pi_{\sf d}$ and $\pi_{\sf 
nd}$ are increasing if and only if $\pi\in{\cal B}_n$ and $\exc(\pi)=\des(\pi)$.
\end{prop}

\begin{bew}
Let $\pi\in{\cal S}_n$ be a permutation for which both 
the descent top word $\pi_{\sf d}$ and the word $\pi_{\sf nd}$ consisting of 
the remaining letters are increasing. Clearly, $\pi$ avoids the pattern $321$. As mentioned 
at the very beginning, hence every descent is an excedance of $\pi$. But these are all 
excedances: assume that $i,i+1,\ldots,i+k\in{\sf E}(\pi)$ and $i-1,i+k+1\notin{\sf 
E}(\pi)$ for some $k>0$. Obviously, there exists an integer $j>i$ (even 
$j>i+k$) satisfying $\pi_j\le i$. Then $j$ can not be a non-descent since $i$ 
is such a one and $\pi_j<\pi_i$. On the other hand, $j$ is not a descent as 
well since $i+k$ is a descent and $\pi_j<\pi_{i+k}$. The converse is evident (${\sf E}(\pi)$ equals ${\sf D}(\pi)$). 
\end{bew}

\begin{rem}
\brm
In \cite[Th. 2]{foata-zeilberger} a natural expression for Denert's 
permutation statistic $\den$ was given: $\den(\pi)=\inv(\pi_{\sf e})+\inv(\pi_{\sf 
ne})+i_1+\ldots+i_e$ where $i_1,\ldots,i_e$ denote the excedances of $\pi\in{\cal 
S}_n$. In case of bi-increasing permutations, $\den$ equals the sum of 
excedances.\\
Since every descent of a bi-increasing permutation is an excedance as well, the conditions $\den(\pi)=\maj(\pi)$ and $\exc(\pi)=\des(\pi)$ are 
equivalent in ${\cal B}_n$. 
\erm
\end{rem}

Calculating the number of permutations $\pi\in{\cal B}_n$ satisfying 
$\exc(\pi)=\des(\pi)$ for some small values $n$ yields the first terms of 
another well-known number sequence, the Motzkin numbers. The {\it Motzkin 
numbers} may be defined by
\bdpm
M_0=1,\quad M_n=M_{n-1}+\sum_{i=0}^{n-2} M_iM_{n-2-i}\quad\mbox{for }n\ge1.
\edpm
One of their numerous combinatorial interpretations is the following one: $M_n$ 
counts the number of lattice paths from the origin to $(n,0)$, with steps 
$[1,1]$ (called {\sf u}p-steps), $[1,-1]$ (called {\sf d}own-steps), and 
$[1,0]$, never going below the $x$-axis. 
Such paths are called {\it Motzkin paths}.\\
A modification are the {\it 2-Motzkin paths} which distinguish two kinds of $[1,0]$ 
steps: {\sf s}olid steps and {\sf b}roken steps. We denote the set of 2-Motzkin 
paths of length $n$ by ${\cal M}_n$. By a simple substitution, a 
2-Motzkin path of length $n$ can be transformed into a Dyck path of length 
$2n+2$. Replace a step {\sf u} by two up-steps, a step {\sf d} by two down 
steps, a step {\sf s} by an up-down combination, a step {\sf b} by a down-up combination and adjoin an 
additional up-step at the beginning and an additional down-step at the end. 
In particular, this bijection shows $|{\cal M}_n|=C_{n+1}$.\\[2ex] 
Together with a weight function, the set ${\cal M}_n$ corresponds to the 
symmetric group ${\cal S}_n$. Several authors (\cite{francon-viennot}, \cite{foata-zeilberger}, \cite{biane}) 
gave bijections regarding this whose connections was studied in \cite{clarke 
etal}.\\
We will use the first one, due to Fran\c{c}on and Viennot, in the notation of 
\cite{clarke etal} to enumerate the bi-increasing permutations for which excedance number 
and descent number are equally.\\[2ex]
Given a permutation $\pi\in{\cal S}_n$ we separate $\pi$ into its {\it descent 
blocks} by putting in a dash between the letters $\pi_i$ and $\pi_{i+1}$ 
whenever $i$ is a non-descent. For example, the permutation 
$\pi=2\:7\:4\:3\:1\:6\:5\:8\:10\:9\in{\cal S}_{10}$ has the descent block decomposition 
$2-7\:4\:3\:1-6\:5-8-10\:9$.\\
Clearly, if $\pi$ is bi-increasing (or $321$-avoiding) then the maximum length of a 
descent block is two.    

\begin{theo} \label{special permutation-path}
Let $\pi\in{\cal B}_n$ satisfy $\exc(\pi)=\des(\pi)=k$. Then $\pi$ corresponds in 
an one-to-one fashion to a Motzkin path of length $n$ with $k$ up-steps.
\end{theo}

\begin{bew}
Given a permutation $\pi\in{\cal S}_n$, Fran\c{c}on-Viennot's bijection defines the 2-Motzkin path 
$c=c_1\cdots c_n$ (which corresponds together with its weight to $\pi$) as follows:
\bdpm
c_{\pi_i}=\left\{\ba{ccl}
{\sf s}&&\mbox{if $\pi_{i-1}<\pi_i<\pi_{i+1}$},\\
{\sf b}&&\mbox{if $\pi_{i-1}>\pi_i>\pi_{i+1}$},\\
{\sf d}&&\mbox{if $\pi_{i-1}<\pi_i>\pi_{i+1}$},\\
{\sf u}&&\mbox{if $\pi_{i-1}>\pi_i<\pi_{i+1}$}
\ea\right.
\edpm
where $\pi_0:=0$ and $\pi_{n+1}:=n+1$. Using the descent block decomposition, 
this means: if the letter $i$ is the first (last) one in a descent block of length at least 
2 then the $i$th step is a down-step (up-step). If $i$ is a letter lain strictly inside a block 
then the $i$th step is broken. If $i$ belongs to a singleton block then the 
$i$th step is solid. (For example, the permutation $2\:7\:4\:3\:1\:6\:5\:8\:10\:9$ 
is associated with the path ${\sf usbbuddsud}$.)\\
Let $\pi\in{\cal B}_n$ with $\exc(\pi)=\des(\pi)=k$.  
Obviously, $c$ has no broken step. Since the words $\pi_{\sf d}$ and $\pi_{\sf 
nd}$ are increasing now the map which takes any bi-increasing permutation to a 
Motzkin path is a bijection. Any descent of $\pi$ corresponds to a down-step. Hence the descent 
number of $\pi$ equals the number of down-steps (or, equivalently, up-steps) in 
$c$. 
\end{bew}
\vspace*{1ex}

Fran\c{c}on-Viennot's map also transforms the descent difference into a 
natural statistic of the 2-Motzkin path. 

\begin{prop} \label{height sum}
Let $\pi\in{\cal S}_n$, and let $c$ be the $2$-Motzkin path to which $\pi$ is taken by 
the map described above. Denote by $h_i$ the {\rm height} of the $i$th step of $c$ defined to be the ordinate of its starting point. 
Then we have $\ddes(\pi)=h_1+\ldots+h_n$.
\end{prop}

\begin{bew}
It is evident that $\ddes(\pi)$ is just the difference of the sum of such letters 
which are the first ones and the sum of those which are the last ones in a block of length at least 2 
obtained by the descent block decomposition of $\pi$. Each of the first 
mentioned letters corresponds to a down-step while each of the last mentioned 
ones corresponds to an up-step of $c$. It is easy to see that 
\bdpm
\sum_{c_i={\sf d}} i-\sum_{c_i={\sf u}} i=h_1+\ldots+h_n.
\edpm
\end{bew}

\begin{exam}
\brm                                        
For $\pi=3\:1\:7\:2\:4\:5\:6\:8\:10\:9\in{\cal B}_{10}$ we have 
$\exc(\pi)=\des(\pi)=3$. Following the proof, we obtain {\sf uudsssdsud} as 
corresponding Motzkin path. Conversely, given the path
\addtocounter{pic}{1}
\begin{center}                                                        
\unitlength0.4cm
\begin{picture}(10,2)
\linethickness{0.5pt}
\bezier{100}(0,0)(1,1)(2,2)
\bezier{50}(2,2)(2.5,1.5)(3,1)
\bezier{150}(3,1)(4.5,1)(6,1)
\bezier{50}(6,1)(6.5,0.5)(7,0)
\bezier{50}(7,0)(7.5,0)(8,0)
\bezier{50}(8,0)(8.5,0.5)(9,1)
\bezier{50}(9,1)(9.5,0.5)(10,0)
\linethickness{0.2pt}
\multiput(0,0)(0,1){3}{\line(1,0){10}}
\multiput(0,0)(1,0){11}{\line(0,1){2}}
\end{picture}

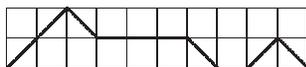
{\footnotesize{\bf Figure \thepic}\hspace*{0.25cm}Motzkin path corresponding to $3\:1\:7\:2\:4\:5\:6\:8\:10\:9$}
\end{center}
we can retrieve $\pi$ as follows: first form blocks consisting of the indices of the 
$k$th down-step and the $k$th up-step. We obtain 
$3\:1-7\:2-10\:9$. Then insert the remaining numbers as singleton blocks such that 
the blocks are increasing ordered by their last letter.
\erm
\end{exam}

A refinement of Motzkin path enumeration according to the number of up-steps 
was done in \cite{donaghey-shapiro}. Using this result, Theorem \ref{special 
permutation-path} yields

\begin{cor}
The number of bi-increasing permutations $\pi\in{\cal S}_B$ satisfying 
$\exc(\pi)=\des(\pi)$ equals $M_n$. In particular, there are ${n\choose 2k}C_k$ 
such permutations with $k$ descents.
\end{cor}

\begin{rems} \label{barcucci and deutsch}
\brm
\begin{enum}
\item[]
\item In \cite{barcucci etal1}, the authors constructed classes of permutations 
of length $n$ which avoid certain patterns. The enumeration of these permutations yields 
an integer sequence whose first term is the $n$th Motzkin number and whose 
limit is the $n$th Catalan number.\\
The first term counts the number of permutations $\pi\in{\cal S}_n$ which avoid 
the patterns $321$ and $3\bar{1}42$. The latter one means that any subsequence 
of type $231$ in $\pi$ must be contained in a subsequence of type $3142$.\\
These are exactly the permutations considered above. Let $i<j<k$ such that 
$\pi_k<\pi_i<\pi_j$ where $\pi\in{\cal B}_n$ with $\exc(\pi)=\des(\pi)$. 
Clearly, $i$ and $j$ are excedances (or, equivalently, descents) but $k$ is none. 
Since $\pi$ avoids $321$ there is a non-excedance $l$ for which $i<l<j$ and 
$\pi_l<\pi_k$. (The second condition follows from Proposition \ref{increasing descent and non-descent 
words}.)\\
Note that the integer sequences $(a_k(n))$ whose $k$th term is defined as number 
of bi-increasing permutations $\pi\in{\cal B}_n$ for which 
$\exc(\pi)-\des(\pi)\le k-1$ are of a similar behaviour as the sequences in \cite{barcucci 
etal1}. We have $a_1(n)=M_n,\;a_{n-2}(n)=C_n-1$, and $a_k(n)=C_n$ for all $k\ge 
n-1$. 
\item The one-to-one correspondence \ref{special permutation-path} yields a simple proof for an observation made by Deutsch 
(\cite{deutsch}). For $n>1$ there are as many Motzkin paths of length $n$ 
with no horizontal steps on the $x$-axis as Motzkin paths of length $n-1$ 
with at least one horizontal step on the $x$-axis.\\
By the construction, the Motzkin path corresponding to $\pi\in{\cal B}_n$ 
connects the lattice points $(i-1,0)$ and $(i,0)$ by a solid step if and only if $i$ is a fixed point of 
$\pi$. (See also Corollary \ref{fixed points}.)\\ 
Let $\sigma\in{\cal B}_{n-1}$ satisfy $\exc(\sigma)=\des(\sigma)$ with 
minimal fixed point $i$. Define the permutation $\pi$ to be obtained from 
$\sigma$ by inserting the letter $n$ between $\sigma_{i-1}$ and $\sigma_{i}$ 
and sorting the descents such that $\pi_{\sf d}$ is increasing. For example, 
for $\sigma=2\:1\:3\:6\:4\:5\:7$ we obtain 
$\pi=\underline{2}\:1\:\underline{6}\:3\:\underline{8}\:4\:5\:7$ (the 
underlined subword is just $\pi_{\sf d}$). Obviously, $\pi\in{\cal B}_n$.
Any fixed point forms a singleton descent block. By inserting $n$ before 
$\sigma_{i}$, a new descent arises. Note that sorting $\pi_{\sf d}$ preserves 
all descents and non-descents since $\sigma_{\sf d}$ and $\sigma_{\sf nd}$ are 
increasing. Moreover, $\pi$ is a derangement: if $\sigma_j=j+1$ for some $j>i$ 
then the letter $\sigma_j$ is moved to the left while sorting since $j$ is a descent 
(equally, excedance).\\
In particular, we have $\exc(\pi)=\exc(\sigma)+1$. Consequently, we can refine Deutsch's statement 
regarding the number of up-steps.
\end{enum}
\erm
\end{rems}

By a little additional convention we can extend the correspondence \ref{special 
permutation-path} to a bijection which takes {\it any} bi-increasing 
permutation to a lattice path, with preserving all parameters as above.

\begin{theo} \label{permutation-path}
There is a bijection between bi-increasing permutations $\pi\in{\cal B}_n$ with 
$\des(\pi)=d$ and $\ddes(\pi)=k$ and $2$-Motzkin paths of length $n$ having $d$ 
up-steps and no broken step of height $0$ such that the sum of heights of all steps equals $k$.
\end{theo}

\begin{bew}
The only difference from the situation studied above is that the letter forming a 
descent block of length one may be an excedance letter now. We encode this 
information by the path as follows. With the convention $\pi_0:=0$ and $\pi_{n+1}:=n+1$ 
we set 
\bdpm
c_{\pi_i}=\left\{\ba{ccl}
{\sf s}&&\mbox{if $\pi_{i-1}<\pi_i<\pi_{i+1}$ and $i$ is not an excedance},\\
{\sf b}&&\mbox{if $\pi_{i-1}<\pi_i<\pi_{i+1}$ and $i$ is an excedance},\\
{\sf d}&&\mbox{if $\pi_{i-1}<\pi_i>\pi_{i+1}$},\\
{\sf u}&&\mbox{if $\pi_{i-1}>\pi_i<\pi_{i+1}$.}
\ea\right.
\edpm
By Remark \ref{barcucci and deutsch}b), the horizontal steps at level zero 
corresponds to fixed points which are clearly non-excedances. The 
transformation of the permutation statistics follows from Theorem \ref{special 
permutation-path} and from the proof of Proposition \ref{height sum}. In particular, the number 
of broken steps equals $\exc(\pi)-\des(\pi)$. 
\end{bew}

This yields another combinatorial interpretation for the Catalan numbers. The 
second statement results from \cite[Cor. 3.3]{robertson etal} (Fine numbers 
count bi-increasing derangements). The generalization in part c) follows from 
Theorem 7.5 in the same paper. Clearly, $\sum_{k=0}^n m_{n,k}=C_n$. 
Since $|{\cal M}_n|=C_{n+1}$, we have in addition $\sum_{k=0}^n 2^km_{n,k}=C_{n+1}$.
\newpage

\begin{cor}
\begin{enum}
\item[]
\item The number of $2$-Motzkin paths of length $n$ having no broken steps on the $x$-axis 
is the $n$th Catalan number $C_n$.
\item The number of $2$-Motzkin paths of length $n$ having no horizontal steps on the $x$-axis 
is the $n$th Fine number $F_n$.
\item The number of $2$-Motzkin paths of length $n$ having no broken steps but 
$k$ solid steps on the $x$-axis equals 
\bdpm
m_{n,k}=\sum_{i=0}^{n-k} (-1)^i\left(\frac{k+1+i}{n+1}\right){2n-k-i\choose 
n}{k+i\choose k}.
\edpm
\end{enum}
\end{cor}

\begin{exam}
\brm
The running example of the previous section $\pi=2\:6\:1\:3\:7\:4\:5\:8\:10\:9\in{\cal B}_{10}$ 
is taken to the 2-Motzkin path
\addtocounter{pic}{1}
\begin{center}                                                        
\unitlength0.4cm
\begin{picture}(10,2)
\linethickness{0.5pt}
\bezier{50}(0,0)(0.5,0.5)(1,1)
\bezier{5}(1,1)(1.5,1)(2,1)
\bezier{50}(2,1)(2.5,1)(3,1)
\bezier{50}(3,1)(3.5,1.5)(4,2)
\bezier{50}(4,2)(4.5,2)(5,2)
\bezier{100}(5,2)(6,1)(7,0)
\bezier{50}(7,0)(7.5,0)(8,0)
\bezier{50}(8,0)(8.5,0.5)(9,1)
\bezier{50}(9,1)(9.5,0.5)(10,0)
\linethickness{0.2pt}
\multiput(0,0)(0,1){3}{\line(1,0){10}}
\multiput(0,0)(1,0){11}{\line(0,1){2}}
\end{picture}

{\footnotesize{\bf Figure \thepic}\hspace*{0.25cm}2-Motzkin path corresponding 
to $2\:6\:1\:3\:7\:4\:5\:8\:10\:9$}
\end{center}
having $\des(\pi)=3$ up-steps and height sum $\ddes(\pi)=9$.
\erm
\end{exam}

Let ${\cal M}^*_n$ be the set of all 2-Motzkin paths of length $n$ whose broken 
steps are all of positive height. The enumeration of bi-increasing permutations 
of length $n$ according to their descent number is equivalent to the enumeration of paths in ${\cal M}^*_n$ 
with regard to the number of the up-steps. The explicit determination of these 
numbers remains open.\\[2ex]
In the previous section we exhibit a one-to-one correspondence between bi-increasing 
permutations and parallelogram polyominoes that transfers $\exc$ and $\dexc$ 
to natural polyomino statistics.\\
The bijection between the symmetric group and weighted 2-Motzkin paths due to Foata and Zeilberger (see \cite{foata-zeilberger}) yields a 
correspondence between ${\cal B}_n$ and ${\cal M}^*_n$ that expresses the 
excedance-based statistics as path parameters, too.\\
Given a permutation $\pi\in{\cal S}_n$, Foata-Zeilberger's map defines the 
2-Motzkin path $c=c_1\cdots c_n$ by
\bdpm
c_i=\left\{\ba{ccl}
{\sf b}&&\mbox{if $i$ is both an excedance and an excedance letter},\\
{\sf s}&&\mbox{if $i$ is both a non-excedance and a non-excedance letter},\\
{\sf u}&&\mbox{if $i$ is both an excedance and a non-excedance letter},\\
{\sf d}&&\mbox{if $i$ is both a non-excedance and an excedance letter}.
\ea\right.
\edpm 
By the construction, no broken step can be of height $0$; hence $c(\pi)\in{\cal 
M}^*_n$ for any $\pi\in{\cal S}_n$. The reduction on bi-increasing permutations yields a bijection 
between ${\cal B}_n$ and ${\cal M}^*_n$. From the indices of up-steps and 
broken steps we obtain the excedances, from the indices of down-steps and 
broken steps we obtain the excedance letters. For a bi-increasing permutation, 
these informations are enough to determine the permutation completely.\\
Obviously, the excedance number of $\pi\in{\cal S}_n$ equals the number of {\sf u} and {\sf 
b} in the corresponding path. Moreover, \cite[proof of Th. 10]{clarke etal} showed that $\dexc(\pi)$ 
is just the sum of heights in terms of path.

\begin{exam}
\brm
The permutation $\pi=2\:6\:1\:3\:7\:4\:5\:8\:10\:9\in{\cal B}_{10}$ is 
transformed into the 2-Motzkin path 
\addtocounter{pic}{1}
\begin{center}                                                        
\unitlength0.4cm
\begin{picture}(10,2)
\linethickness{0.5pt}
\bezier{50}(0,0)(0.5,0.5)(1,1)
\bezier{5}(1,1)(1.5,1)(2,1)
\bezier{100}(2,1)(3,1)(4,1)
\bezier{50}(4,1)(4.5,1.5)(5,2)
\bezier{100}(5,2)(6,1)(7,0)
\bezier{50}(7,0)(7.5,0)(8,0)
\bezier{50}(8,0)(8.5,0.5)(9,1)
\bezier{50}(9,1)(9.5,0.5)(10,0)
\linethickness{0.2pt}
\multiput(0,0)(0,1){3}{\line(1,0){10}}
\multiput(0,0)(1,0){11}{\line(0,1){2}}
\end{picture}

{\footnotesize{\bf Figure \thepic}\hspace*{0.25cm}2-Motzkin path corresponding 
to $2\:6\:1\:3\:7\:4\:5\:8\:10\:9$ by Foata-Zeilberger}
\end{center}
whose heights amount to $\dexc(\pi)=8$ and which contains $\exc(\pi)=4$ steps 
{\sf u} or {\sf b}.
\erm
\end{exam}

Using both, the correspondence due to Foata-Zeilberger and the correspondence from 
Theorem \ref{permutation-path}, shows trivially the equidistribution of 
$\dexc$ and $\ddes$ on the set ${\cal B}_n$. In the next section we will give a 
simple combinatorial proof in terms of bi-increasing permutations for this fact.
\vspace*{0.75cm}


\setcounter{section}{5}\setcounter{theo}{0}

\centerline{\large{\bf 5}\hspace*{0.25cm}
{\sc Parallelogram polyominoes and 2-Motzkin paths}}
\vspace*{0.5cm}

In contrast to the excedance number and descent number, the difference statistics are also equidistributed 
over the set of bi-increasing permutations. But the bijection $\phi$ from 
Proposition \ref{foata} does not furnish proof because ${\cal B}_n$ is not closed under $\phi$.

\begin{theo}
We have $B_n^{\dexc}(k)=B_n^{\ddes}(k)$ for all $n$ and $k$.
\end{theo}

\begin{bew}
We give a bijection $\psi:{\cal B}_n\to{\cal B}_n$ which takes a 
bi-increasing permutation $\pi$ to a bi-increasing permutation $\sigma$ whose excedance letters 
equal those of $\pi$.\\
Given $\pi\in{\cal B}_n$, let $i_1,\ldots,i_s$ be the excedances which are 
descents in addition, and let $j_1,\ldots,j_t$ be the remaining excedances of 
$\pi$. (As mentioned above, $\des(\pi)=s$.) The permutation $\sigma$ is defined as the bi-increasing one whose 
excedance set is $\{\pi_{i_1+1},\ldots,\pi_{i_s+1},\pi_{j_1},\ldots,\pi_{j_t}\}$ 
and whose excedance letters are $\pi_{i_1},\ldots,\pi_{i_s},\pi_{j_1},\ldots,\pi_{j_t}$. 
Note that $\sigma$ is well defined. To see this, let $j$ be the first one of a 
sequence of consecutive excedances whose last term, say $i$, is a descent in addition. 
Then $\pi_j>\pi_{i+1}$ since there is at least one integer $k>j$ with 
$\pi_k<\pi_j$.
\\Obviously, the map $\psi$ is bijective, and we have
\bdpm
\ddes(\pi)=\sum_{k=1}^s (\pi_{i_k}-\pi_{i_k+1})=
\Big(\sum_{k=1}^s \pi_{i_k}+\sum_{k=1}^t \pi_{j_k}\Big)-
\Big(\sum_{k=1}^s \pi_{i_k+1}+\sum_{k=1}^t \pi_{j_k}\Big)=\dexc(\sigma).
\edpm 
\end{bew}

\begin{exam} \label{prevexam}
\brm
Let $\pi=\overline{2}\:\underline{6}\:1\:3\:\underline{7}\:4\:5\:8\:\underline{10}\:9\in{\cal 
B}_{10}$ again. (The letters that are as well excedance letter as descent top are underlined; the 
letters belonging to only-excedances are overlined.) We obtain the permutation 
$\psi(\pi)=2\:6\:1\:7\:3\:4\:5\:8\:10\:9$.
\erm
\end{exam}

\begin{rems} \label{properties of psi}
\brm
\begin{enum}
\item[]
\item The bijection $\psi$ preserves both the excedance number as the number of 
fixed points. The latter one is an immediate consequence of Corollary 
\ref{fixed points}.
\item Let $c\in{\cal M}_n^*$ be the 2-Motzkin path to which $\pi\in{\cal B}_n$ 
is taken by the bijection given in Theorem \ref{permutation-path}. It is not 
difficult to see that $c$ is precisely the path obtained from 
$\psi(\pi)\in{\cal B}_n$ by Foata-Zeilberger's bijection. 
\end{enum}
\erm
\end{rems}

The $q$-series whose coefficient of $x^nq^k$ is just $B_n^{\ddes}(k)$ is 
given in Theorem \ref{series}. From the combinatorial interpretations of the 
distribution coefficients $B_n^{\dexc}(k)$ and $B_n^{\ddes}(k)$, respectively, we obtain

\begin{cor}
There are as many parallelogram polyominoes of perimeter $2n+2$ whose 
area minus height equals $k$ as $2$-Motzkin paths in ${\cal M}^*_n$ whose height sum equals $k$.
\end{cor}

Translated into the languages of polyominoes and lattice paths the bijection $\psi$ 
can be read as follows. Let $(\gamma,\delta)$ be a parallelogram polyomino of perimeter $2n+2$ and 
width $w$. (Hence its height equals $n+1-w$.) Define the integer sets 
$A=\{\gamma_1+\ldots+\gamma_i+i:i=1,2,\ldots,w-1\}$ and $B=\{\delta_1+\ldots+\delta_i+i:i=1,2,\ldots,w-1\}$,
and let $C$ be their intersection. The 2-Motzkin path $c\in{\cal M}^*_n$ 
associated with $(\gamma,\delta)$ we obtain by the convention
\bdpm
c_i=\left\{\ba{ccl}
{\sf d}&&\mbox{if $i\in A\setminus C$},\\
{\sf u}&&\mbox{if $i\in B\setminus C$},\\
{\sf b}&&\mbox{if $i\in C$},\\
{\sf s}&&\mbox{if $i\notin A\cup B$.}
\ea\right.
\edpm
The following correspondences under this bijection are clear from the theorems 
\ref{permutation-parallelogram polyomino} and \ref{permutation-path} and Remark 
\ref{properties of psi}a:
\newpage
\beas
\mbox{\bf parallelogram polyominos}&&\mbox{\bf 2-Motzkin paths}\\
\mbox{perimeter }2n+2&&\mbox{length }n\\
\mbox{width}&&1+\mbox{number of {\sf u}'s}+\mbox{number of {\sf b}'s}\\
\mbox{(height)}&&(\mbox{number of {\sf u}'s}+\mbox{number of {\sf s}'s})\\
\mbox{area minus height}&&\mbox{sum of heights}\\
\mbox{number of rows consisting of one cell}&&\mbox{number of solid steps at level 
$0$}    
\eeas

\begin{exam}
\brm
Let $\gamma=(1,3,0,2,0)$ and $\beta=(0,0,1,4,1)$. Note that $(\gamma,\delta)$  
is just the parallelogram polyomino corresponding to $2\:6\:1\:7\:3\:4\:5\:8\:10\:9\in{\cal 
B}_{10}$, the permutation $\psi(\pi)$ from the previous example. We have 
$A=\{2,6,7,10\},\;B=\{1,2,4,9\}$, and $C=\{2\}$. Thus we obtain the following 
path: 
\addtocounter{pic}{1}
\begin{center}                                                        
\unitlength0.4cm
\begin{picture}(20,6)
\linethickness{0.5pt}
\put(0,0){\line(1,0){3}}
\put(0,1){\line(1,0){4}}
\multiput(1,2)(0,1){3}{\line(1,0){3}}
\multiput(3,5)(0,1){2}{\line(1,0){2}}
\put(0,0){\line(0,1){1}}
\multiput(1,0)(1,0){2}{\line(0,1){4}}
\put(3,0){\line(0,1){6}}
\put(4,1){\line(0,1){5}}
\put(5,5){\line(0,1){1}}
\bezier{50}(10,2)(10.5,2.5)(11,3)
\bezier{5}(11,3)(11.5,3)(12,3)
\bezier{50}(12,3)(12.5,3)(13,3)
\bezier{50}(13,3)(13.5,3.5)(14,4)
\bezier{50}(14,4)(14.5,4)(15,4)
\bezier{100}(15,4)(16,3)(17,2)
\bezier{50}(17,2)(17.5,2)(18,2)
\bezier{50}(18,2)(18.5,2.5)(19,3)
\bezier{50}(19,3)(19.5,2.5)(20,2)
\linethickness{0.2pt}
\multiput(10,2)(0,1){3}{\line(1,0){10}}
\multiput(10,2)(1,0){11}{\line(0,1){2}}
\end{picture}

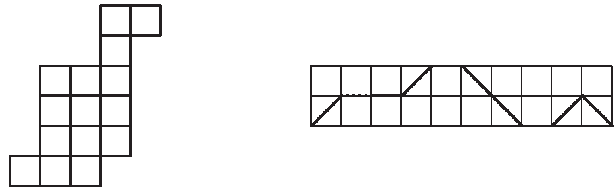
{\footnotesize{\bf Figure \thepic}\hspace*{0.25cm}Parallelogram polyomino and 
corresponding 2-Motzkin path}
\end{center}
(Recall that this is exactly the 2-Motzkin path associated with the permutation 
$\pi\in{\cal B}_{10}$ from Example \ref{prevexam}.)
\erm
\end{exam}

Parallelogram polyominoes can be interpreted as connected skew diagrams. For 
two partitions $\lambda$ and $\mu$ with $\mu\subseteq\lambda$ (i.e., $\mu_i\le\lambda_i$ for all 
$i$), the {\it skew diagram} $\lambda/\mu$ is defined to be the set theoretical difference of the Young 
diagrams associated with $\lambda$ and $\mu$, respectively. The pair $(\gamma,\delta)$ 
describes the skew diagram $1^{\delta_1}2^{\delta_2}\cdots w^{\delta_w}/1^{\gamma_2}2^{\gamma_3}\cdots 
(w-1)^{\gamma_w}$.\\ 
In \cite{nazarov-tarasov}, the rank of a skew diagram was 
introduced. For a skew diagram $\lambda/\mu$, its {\it rank} is defined as difference 
of outside diagonal lengths and inside diagonal lengths and is denoted by $\rk(\lambda/\mu)$.\\
For example, the skew diagram $3^24^35/1^33^2$ (or $((1,3,0,2,0),(0,0,2,3,1))$ in the previous notation)
\addtocounter{pic}{1}
\begin{center}
\unitlength0.4cm
\begin{picture}(5,6)
\linethickness{0.5pt}
\multiput(0,0)(0,1){2}{\line(1,0){3}}
\multiput(1,2)(0,1){3}{\line(1,0){3}}
\multiput(3,5)(0,1){2}{\line(1,0){2}}
\put(0,0){\line(0,1){1}}
\multiput(1,0)(1,0){2}{\line(0,1){4}}
\put(3,0){\line(0,1){6}}
\put(4,2){\line(0,1){4}}
\put(5,5){\line(0,1){1}}
\put(3.5,5.5){\makebox(0,0)[cc]{\footnotesize$+$}}
\multiput(1.5,3.5)(1,-1){2}{\makebox(0,0)[cc]{\footnotesize$+$}}
\put(0.5,0.5){\makebox(0,0)[cc]{\footnotesize$+$}}
\put(1.5,0.5){\makebox(0,0)[cc]{\footnotesize$-$}}
\put(3.5,3.5){\makebox(0,0)[cc]{\footnotesize$-$}}
\end{picture}

{\footnotesize{\bf Figure \thepic}\hspace*{0.25cm}A skew diagram and its rank}
\end{center}
has rank $2$. (All squares belonging to an outside diagonal are marked by $+$, the inside diagonal squares by $-$.)\\ 
If $\mu=\emptyset$, i.e., if $\lambda/\mu$ is a partition, then $\lambda/\mu$ has 
only one outside diagonal (the main diagonal) and no inside diagonals. Hence in this case, $\rk(\lambda/\mu)$ is just 
the {\it Durfee rank} defined for partitions. 
Stanley gave in \cite{stanley2} several equivalent definitions of 
$\rk(\lambda/\mu)$. 
\begin{prop}
Let $(\gamma,\delta)$ be a parallelogram polyomino and $c\in{\cal M}^*_n$ its corresponding 
$2$-Motzkin path. Then the rank of $(\gamma,\delta)$ equals 
\bdpm
1+\mbox{\rm number of {\sf dd}}+\mbox{\rm number of {\sf db}}+
\mbox{\rm number of {\sf sd}}+\mbox{\rm number of {\sf sb}}
\edpm
appearing in $c$.
\end{prop}

\begin{bew}
We use the rank definition dealing with the reduced code of skew diagrams. 
Given a skew diagram $\lambda/\mu$ of perimeter $2n+2$, we mark every vertical boundary 
part by $0$ and every horizontal boundary part by $1$. Reading these numbers 
while moving north and east along the lower boundary (starting from the 
left-hand edge) we obtain a binary sequence $a(\lambda/\mu)=a_1a_2\cdots a_{n+1}$. 
In a similar way, when we read the labels as we move north and east along the upper 
boundary we obtain a binary sequence $b(\lambda/\mu)=b_1b_2\cdots b_{n+1}$. The 
two-line array 
\bdpm
{\rm cd}(\lambda/\mu):=\ba{cccc}a_1&a_2&\cdots&a_{n+1}\\[-1ex]b_1&b_2&\cdots&b_{n+1}\ea
\edpm
we call the {\it reduced code} of $\lambda/\mu$. (Clearly, $a_1=b_{n+1}=1$ and 
$a_{n+1}=b_1=0$.) The concept should be clear from the example:
\addtocounter{pic}{1}
\begin{center}
\unitlength0.4cm
\begin{picture}(21,6.2)
\linethickness{0.3pt}
\multiput(0,0)(0,1){2}{\line(1,0){3}}
\multiput(1,2)(0,1){3}{\line(1,0){3}}
\multiput(3,5)(0,1){2}{\line(1,0){2}}
\put(0,0){\line(0,1){1}}
\multiput(1,0)(1,0){2}{\line(0,1){4}}
\put(3,0){\line(0,1){6}}
\put(4,2){\line(0,1){4}}
\put(5,5){\line(0,1){1}}
\multiput(3.5,6)(1,0){2}{\makebox(0,0)[cc]{\bf\tiny 1}}
\multiput(1.5,4)(1,0){2}{\makebox(0,0)[cc]{\bf\tiny 1}}
\put(0.5,1){\makebox(0,0)[cc]{\bf\tiny 1}}
\put(4.5,5){\makebox(0,0)[cc]{\bf\tiny 1}}
\put(3.5,2){\makebox(0,0)[cc]{\bf\tiny 1}}
\multiput(0.5,0)(1,0){3}{\makebox(0,0)[cc]{\bf\tiny 1}}
\put(0,0.5){\makebox(0,0)[cc]{\bf\tiny 0}}
\put(5,5.5){\makebox(0,0)[cc]{\bf\tiny 0}}
\multiput(1,1.5)(0,1){3}{\makebox(0,0)[cc]{\bf\tiny 0}}
\multiput(3,4.5)(0,1){2}{\makebox(0,0)[cc]{\bf\tiny 0}}
\multiput(4,2.5)(0,1){3}{\makebox(0,0)[cc]{\bf\tiny 0}}
\multiput(3,0.5)(0,1){2}{\makebox(0,0)[cc]{\bf\tiny 0}}
\put(7,3){\makebox(0,0)[lc]{${\rm cd}(3^24^35/1^33^2)=
{1\atop0}{1\atop1}{1\atop0}{0\atop0}{0\atop0}
{1\atop1}{0\atop1}{0\atop0}{0\atop0}{1\atop1}{0\atop1}$}}
\end{picture}

{\footnotesize{\bf Figure \thepic}\hspace*{0.25cm}Reduced code of a skew diagram}
\end{center}
By \cite[Prop. 2.2]{stanley2}, $\rk(\lambda/\mu)$ equals the number of columns 
${0\atop1}$ (or, equivalently, ${1\atop0}$) of ${\rm cd}(\lambda/\mu)$. (The figure shows a skew diagram of rank $2$.)\\
Consider now the reduced code of $(\gamma,\delta)$. In the binary sequence $b$ 
the $1$'s appear at the positions $\gamma_1+1,\gamma_1+\gamma_2+2,\ldots,
\gamma_1+\ldots+\gamma_{w-1}+w-1,\gamma_1+\ldots+\gamma_w+w=n+1$. These 
integers, excepting $n+1$, are exactly the elements of the set $A$ appearing in 
the definition of $c$. The sequence $a$ contains a $1$ at each of the positions 
$1,\delta_1+2,\delta_1+\delta_2+3,\ldots,\delta_1+\ldots+\delta_{w-1}+w$. 
Apart from $1$, decreasing these integers by 1 yields the elements of 
$B$. Thus the $i$th column of the code equals ${0\atop1}$ if and only if 
$i\in A$ and $i-1\notin B$ or $i=n+1$. The first condition means the occurrence 
of {\sf dd}, {\sf sd}, {\sf db} and {\sf sb} in the path $c$.
\end{bew}

\begin{rem} \label{symmetry of rank}
\brm
Because of the symmetry, the step combinations ${\sf uu},{\sf us},{\sf bu}$ and 
${\sf bs}$ (corresponding to code columns ${1\atop0}$) can be counted as well. 
Note that the reverse 2-Motzkin path $c_nc_{n-1}\cdots c_1$ corresponds to the 
parallelogram polyomino obtained by rotating $(\gamma,\delta)$ $180^{\circ}$. 
\erm
\end{rem}

We will apply the bijection to enumerate partitions of prescribed perimeter according their rank. 
(Here the perimeter of a partition is just twice the sum of its largest part 
and the number of its parts.)

\begin{lem} \label{Motzkin path for partitions}
Any $2$-Motzkin path $c\in{\cal M}^*_n$ corresponds to a partition if and only if 
there exists an integer $k\in[n]$ such that $c_1,\ldots,c_k\in\{{\sf u},{\sf s}\}$ and 
$c_{k+1},\ldots,c_n\in\{{\sf d},{\sf b}\}$.
\end{lem}

\begin{bew} 
Let $(\gamma,\delta)$ be a parallelogram polyomino of perimeter $2n+2$ and 
width $w$. (Recall that 
$\gamma_1+\ldots+\gamma_w=\delta_1+\ldots+\delta_w=n+1-w$.) Then 
$(\gamma,\delta)$ describes a partition if and only if 
$\gamma_2=\ldots=\gamma_w=0$. Thus the elements of the set $A$ appearing in 
the bijection definition are exactly the integers $n+2-w,n+3-w,\ldots,n$ . Consequently, $c_i\in\{{\sf d},{\sf 
b}\}$ if $i\ge n+2-w$ and $c_i\in\{{\sf u},{\sf s}\}$ otherwise.
\end{bew}

\begin{lem} \label{binomial-identity}
For all non-negative integers $a<b\le c$ we have
\bdpm
\sum_{k=0}^{c-b} {a+k\choose a}{c-a-1-k\choose b-a-1}={c\choose b}.  
\edpm
\end{lem}

\begin{bew}
This special case of the Chu-Vandermonde identity can be proofed very easily. 
Draw a line of $c$ dots and circle $b$ dots. Then count the number of uncircled 
dots to the left of the $(a+1)$st circled dot.
\end{bew}

\begin{theo}
There are ${n\choose 2r-1}$ partitions of perimeter $2n+2$ and rank $r$.
\end{theo}

\begin{bew}
Let $\lambda$ be a partition of perimeter $2n+2$ and rank $r$, and $c_1\cdots 
c_n$ the corresponding $2$-Motzkin path. By Lemma \ref{Motzkin path for partitions} and Remark \ref{symmetry of rank}, 
the rank is one plus the number of occurrences of {\sf uu} and {\sf us} in $c_1\cdots 
c_k$ where $k$ is the maximum integer for which $c_i\in\{\sf u,s\}$. Clearly, there are 
${k-1\choose r-1}$ different words $c_1\cdots c_k$ over $\{\sf u,s\}$ containing 
$r-1$ times {\sf uu} or {\sf us}. Obviously, $c_{k+1}\cdots c_n$ has $r-1$ letters {\sf d} (as 
counterparts to $r-1$ {\sf u}'s in the first part of $c$), and only {\sf b} 
otherwise. There are ${n-k\choose r-1}$ such words. Consequently, the desired 
number equals
\bdpm
\sum_{k=r}^{n+1-r} {k-1\choose r-1}{n-k\choose r-1}=
\sum_{k=0}^{n+1-2r} {k+r-1\choose r-1}{n-k-r\choose r-1}={n\choose 2r-1}.
\edpm
(For the second identity use Lemma \ref{binomial-identity} taking $a=r-1,\;b=2r-1,\;c=n$.)
\end{bew}

To prove the generalization we utilize an encoding principle introduced in 
\cite[Sect. 5]{deutsch-shapiro}. Given a parallelogram polyomino of perimeter 
$2n+2$, the corresponding 2-Motzkin path $c\in{\cal M}_{n+1}$ is defined as 
follows. Beginning at the origin read the upper border path and lower border path, respectively, 
step by step. Reading the $i$th steps, set  
\bdpm
c_i=\left\{\ba{ccl}
{\sf u}&&\mbox{if upper path goes north, lower path goes east},\\
{\sf d}&&\mbox{if upper path goes east, lower path goes north},\\
{\sf b}&&\mbox{if both paths go east},\\
{\sf s}&&\mbox{if both paths go north.}
\ea\right.
\edpm
For instance, the parallelogram polyomino 
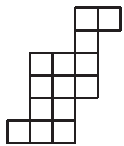
\begin{center}
\unitlength0.3cm
\begin{picture}(5,6)
\linethickness{0.5pt}
\multiput(0,0)(0,1){2}{\line(1,0){3}}
\multiput(1,2)(0,1){3}{\line(1,0){3}}
\multiput(3,5)(0,1){2}{\line(1,0){2}}
\put(0,0){\line(0,1){1}}
\multiput(1,0)(1,0){2}{\line(0,1){4}}
\put(3,0){\line(0,1){6}}
\put(4,2){\line(0,1){4}}
\put(5,5){\line(0,1){1}}
\end{picture}
\end{center}
is associated with the 2-Motzkin path 
{\sf ubussbdssbd}. Note that every 2-Motzkin path arising in this way begins with an up-step and ends with a 
down-step. Moreover, all steps are of height at least one, excepting the first 
one. Thus by deleting the very first and very last step of $c$ we obtain a 
correspondence between parallelogram polyominoes of perimeter $2n+2$ and 
2-Motzkin path of length $n-1$ which is one-to-one. 

\begin{theo}
The number of all connected skew diagrams of perimeter $2n+2$ and rank $r$ 
equals to 
\bdpm
2^{n+1-2r}{n-1\choose 2r-2}C_{r-1}.
\edpm
\end{theo}

\begin{bew}
Let $\lambda/\mu$ be a connected skew diagram of perimeter $2n+2$ and $c$ its 
corresponding 2-Motzkin path of length $n+1$ obtained as described above.
By construction, the $i$th column of the reduced code of $\lambda/\mu$ equals ${0\atop1}$ if and only if 
the $i$th step of the lower border path goes north, and the $i$th step of the upper 
border path goes east. Hence each such column corresponds to a down-step of $c$.
Consequently, Deutsch-Shapiro's bijection takes a skew diagram of perimeter 
$2n+2$ and rank $r$ to a 2-Motzkin path of length $n-1$ having $r-1$ 
down-steps. (The last down-step of $c$ has been deleted.) By 
\cite{donaghey-shapiro}, the number of Motzkin paths of length $n-1$ with $r-1$ down-steps equals ${n-1\choose 
2r-2}C_{r-1}$. Thus the number given in 
the statement counts the 2-Motzkin paths with these parameters.
\end{bew}

\begin{cor}
There are as many $2$-Motzkin paths in ${\cal M}_{n-1}$ with $k$ down-steps as 
$2$-Motzkin paths in ${\cal M}^*_n$ with $k$ occurrences of double steps {\sf 
dd}, {\sf db}, {\sf sd}, and {\sf sb}.
\end{cor}

\begin{rem}
\brm
The statement of Lemma \ref{Motzkin path for partitions} holds analogously for the 
bijection of Deutsch and Shapiro as well.
\erm
\end{rem}
\vspace*{0.6cm}


\setcounter{section}{6}\setcounter{theo}{0}

\centerline{\large{\bf 6}\hspace*{0.25cm}
{\sc An application to the distributions over the symmetric group}}
\vspace*{0.5cm}

In the first section we introduce the surjection which takes any permutation 
$\pi\in{\cal S}_n$ to the bi-increasing permutation $\hat{\pi}\in{\cal B}_n$. 
The map induces a partition of ${\cal S}_n$ into $C_n$ disjoint classes, each represented 
by a bi-increasing permutation. All the elements of a class have the same 
excedance number and excedance difference. This section deals first with the question: for $\pi\in{\cal 
B}_n$, how many permutations belong to the class $[\pi]$?\\[2ex] 
As discussed in the proof of Corollary \ref{bounds}, every member of $[\pi]$ 
arises from $\pi$ by applying some transpositions to $\pi$ which preserve all 
excedance and their letters. Let $T_\pi$ denote the set of pairs $(i,j)$ such that 
either $i<j$ are excedances and $j<\pi_i$ or $i<j$ are non-excedances and 
$i\ge\pi_j$. Clearly, exchanging $\pi_i$ and $\pi_j$ has no influence on $\exc$ 
and $\dexc$ if and only if $(i,j)\in T_\pi$. In the proof of \ref{bounds} it 
was shown that $|T_\pi|=\dexc(\pi)-\exc(\pi)$.\\
Note that if $(i,j)$ and $(i,k)$ belong to $T_\pi$ then $(j,k)$ does it as well, provided that $j<k$. 
Every permutation of $[\pi]$ can be constructed from $\pi$ by applying a sequence 
of transpositions $(i_1,j_1),\ldots,(i_s,j_s)\in T_\pi$ to the positions of $\pi$ where 
$i_1<i_2<\ldots<i_s$ for some $s\ge 0$. Consequently, we have
\bdpm
|[\pi]|=\prod_{i=1}^{n-1}\Big(|\{j:(i,j)\in T_\pi\}|+1\Big).
\edpm
For example, for $\pi=2\:6\:1\:3\:7\:4\:5\:8\:10\:9\in{\cal B}_{10}$ we have 
$T_\pi=\{(2,5),(3,4),(4,6),(6,7)\}$, and hence $|[\pi]|=16$. In detail, the 
permutations 
\bdpm
\ba{ccccccc}
2\:6\:1\:3\:7\:4\:5\:8\:10\:9&&2\:6\:1\:3\:7\:5\:4\:8\:10\:9&&
2\:6\:3\:4\:7\:1\:5\:8\:10\:9&&2\:7\:3\:1\:6\:5\:4\:8\:10\:9\\[-1ex]
2\:7\:1\:3\:6\:4\:5\:8\:10\:9&&2\:7\:3\:1\:6\:4\:5\:8\:10\:9&&
2\:6\:3\:1\:7\:5\:4\:8\:10\:9&&2\:7\:1\:4\:6\:5\:3\:8\:10\:9\\[-1ex]
2\:6\:3\:1\:7\:4\:5\:8\:10\:9&&2\:7\:1\:4\:6\:3\:5\:8\:10\:9&&
2\:6\:1\:4\:7\:5\:3\:8\:10\:9&&2\:6\:3\:4\:7\:5\:1\:8\:10\:9\\[-1ex]
2\:6\:1\:4\:7\:3\:5\:8\:10\:9&&2\:7\:1\:3\:6\:5\:4\:8\:10\:9&&
2\:7\:3\:4\:6\:1\:5\:8\:10\:9&&2\:7\:3\:4\:6\:5\:1\:8\:10\:9
\ea
\edpm
are precisely the elements of $[\pi]$.\\[2ex]
Surprisingly, the number $|[\pi]|$ can immediately read off from the polyomino 
connected with $\pi\in{\cal B}_n$ by the correspondence given in Section 3. 

\begin{theo}
Let $(\alpha,\beta)\in{\cal C}_{n,w}\times {\cal C}_{n,w}$ be the step polyomino 
corresponding to $\pi\in{\cal B}_n$. Denote by $R_1,\ldots,R_{n+1}$ the rows of 
$(\alpha,\beta)$, and let $a_i$ be the number of columns in common to $R_i$ and 
$R_{i+1}$, for $1\le i\le n$. Then $|[\pi]|=a_1a_2\cdots a_n$.
\end{theo}

\begin{bew}
First we show that 
$a_i=|\{j\in[w-1]:\beta_1+\ldots+\beta_j<i\le\alpha_1+\ldots+\alpha_j\}|+1$.\\
By the bijection given in \ref{permutation-polyomino}, the $j$th column begins 
at level $\beta_1+\ldots+\beta_{j-1}$ where $\beta_0:=0$ and ends at level 
$\alpha_1+\ldots+\alpha_j+1$. Hence both $R_i$ and $R_{i+1}$ contain squares of 
the $j$th column if 
$\beta_1+\ldots+\beta_{j-1}<i<\alpha_1+\ldots+\alpha_j+1$. Thus, $a_i$ is the 
number of $j\in[w]$ satisfying either 
$\beta_1+\ldots+\beta_j<i\le\alpha_1+\ldots+\alpha_j$ or 
$\beta_1+\ldots+\beta_{j-1}<i\le\beta_1+\ldots+\beta_j$. Clearly, for each $i$ 
there exists exactly one integer $j$ for which the second condition holds. 
(Since $\beta_1+\ldots+\beta_w=n=\alpha_1+\ldots+\alpha_w$ we have $j\in[w-1]$.)
\\
The transpositions $(i,\cdot)\in T_\pi$ have already been counted in the proof of Corollary 
\ref{bounds}. Let $b_i$ denote this number, increased by $1$. For any excedance $i$, we 
obtain $b_i=|\{j\in{\sf E}(\pi):i<j<\pi_i\}|+1=|\{j\in{\sf E}(\pi):i\le j<\pi_i\}|$. It is easy to see that
\bdpm
\prod_{i\in{\sf E(\pi)}} |\{j\in{\sf E}(\pi):i\le j<\pi_i\}|
=\prod_{i\in{\sf E(\pi)}} |\{j\in{\sf E}(\pi):j\le i<\pi_j\}|.
\edpm 
On the other hand, if $i$ is a non-excedance then $b_i=|\{j\notin{\sf E}(\pi):\pi_j\le 
i<j\}|+1$, and by Lemma \ref{auxiliary}b, we have  
$b_i=|\{j\in{\sf E}(\pi):j<i<\pi_j\}|+1$.\\
By Remark \ref{partial sums}, the partial sums 
$\beta_1,\beta_1+\beta_2,\ldots,\beta_1+\ldots+\beta_{w-1}$ are just the 
excedances of $\pi$, and $\alpha_1+1,\alpha_1+\alpha_2+1,\ldots,\alpha_1+\ldots+\alpha_{w-1}+1$ their 
letters. Consequently, we obtain 
$b_i=|\{j\in[w-1]:\beta_1+\ldots+\beta_j<i<\alpha_1+\ldots+\alpha_j+1\}|+1$ for 
any non-excedance $i$, and
\bdpm
\prod_{i\in{\sf E(\pi)}} b_i=\prod_{i\in{\sf E(\pi)}} 
(|\{j\in[w-1]:\beta_1+\ldots+\beta_j<i<\alpha_1+\ldots+\alpha_j+1\}|+1)=
\prod_{i\in{\sf E(\pi)}} a_i.
\edpm
\end{bew}

\begin{cor}
Let $(\gamma,\delta)$ be the parallelogram polyomino corresponding to $\pi\in{\cal B}_n$. Denote by 
$a_1,\ldots,a_n$ the diagonal lengths of $(\gamma,\delta)$. Then $|[\pi]|=a_1a_2\cdots a_n$.
\end{cor}

\begin{bew}
Figuratively, the polyomino transformation from Section 3 works as 
follows. Given a step polyomino, move the squares contained in the $k$th column toward the 
bottom, each by $k-1$ units. Then remove the top squares of all columns. Clearly, 
common borders of adjacent rows correspond to the diagonals of the resulting 
polyomino. 
\end{bew}

\begin{exam}
\brm
For $\pi=2\:6\:1\:3\:7\:4\:5\:8\:10\:9\in{\cal B}_{10}$, we have $|[\pi]|=16$
(see above). Alternatively, this can be seen from the corresponding step polyomino 
and parallelogram polyomino, respectively: 
\addtocounter{pic}{1}
\begin{center}                                                        
\unitlength0.4cm
\begin{picture}(22,12)
\linethickness{0.3pt}
\put(1,0){\line(1,0){1}}
\put(1,1){\line(1,0){2}}
\put(1,2){\line(1,0){3}}
\multiput(2,3)(0,1){2}{\line(1,0){2}}
\multiput(2,5)(0,1){2}{\line(1,0){3}}
\put(3,7){\line(1,0){2}}
\put(4,8){\line(1,0){1}}
\multiput(4,9)(0,1){2}{\line(1,0){2}}
\put(5,11){\line(1,0){1}}
\put(1,0){\line(0,1){2}}
\multiput(2,0)(1,1){2}{\line(0,1){6}}
\put(4,2){\line(0,1){8}}
\put(5,5){\line(0,1){6}}
\put(6,9){\line(0,1){2}}
\multiput(13,0)(0,1){2}{\line(1,0){3}}
\multiput(14,2)(0,1){3}{\line(1,0){3}}
\multiput(16,5)(0,1){2}{\line(1,0){2}}
\put(13,0){\line(0,1){1}}
\multiput(14,0)(1,0){2}{\line(0,1){4}}
\put(16,0){\line(0,1){6}}
\put(17,2){\line(0,1){4}}
\put(18,5){\line(0,1){1}}
\multiput(13.5,0.5)(1,0){3}{\circle*{0.3}}
\multiput(14.5,1.5)(1,0){2}{\circle*{0.3}}
\multiput(14.5,2.5)(1,0){3}{\circle*{0.3}}
\multiput(14.5,3.5)(1,0){3}{\circle*{0.3}}
\put(16.5,4.5){\circle*{0.3}}
\multiput(16.5,5.5)(1,0){2}{\circle*{0.3}}
\linethickness{1pt}
\put(1,1){\line(1,0){1}}
\put(2,2){\line(1,0){1}}
\multiput(2,3)(0,1){3}{\line(1,0){2}}
\put(3,6){\line(1,0){2}}
\multiput(4,7)(0,1){3}{\line(1,0){1}}
\put(5,10){\line(1,0){1}}
\linethickness{0.5pt}
\bezier{50}(14.5,1.5)(15,1)(15.5,0.5)
\bezier{50}(14.5,2.5)(15,2)(15.5,1.5)
\bezier{50}(14.5,3.5)(15,3)(15.5,2.5)
\bezier{50}(15.5,3.5)(16,3)(16.5,2.5)
\end{picture}

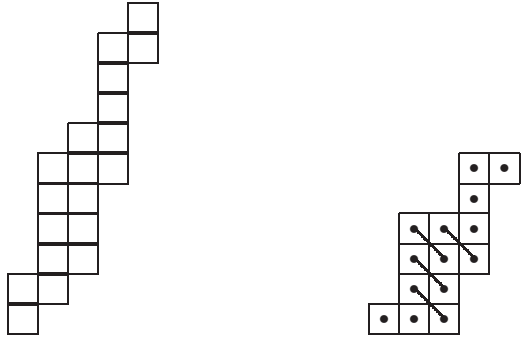
{\footnotesize{\bf Figure \thepic}\hspace*{0.25cm}Size of 
$[2\:6\:1\:3\:7\:4\:5\:8\:10\:9]$ from the polyominoes}
\end{center}
The product of the lengths drawn bold in the left-hand polyomino equals just 
$|[\pi]|$, just as the product of the diagonal lengths for the right-hand 
polyomino.  
\erm
\end{exam}

Corollary \ref{symmetry} said that there are as many bi-increasing permutations 
of length $n$ with $e$ excedances and excedance difference $k$ as those 
having $n-1-e$ excedances and excedance difference $n-1-2e+k$. 
The same result holds when we consider arbitrary permutations.

\begin{cor}
$S_n^{(\exc,\dexc)}(e,k)=S_n^{(\exc,\dexc)}(n-1-e,n-1-2e+k)$ for all $e$ and $k$.
\end{cor}

\begin{bew}
Corollary \ref{symmetry} has been proved by an involution in terms of 
polyominoes. Given $\pi\in{\cal B}_n$ with $\exc(\pi)=e$ and $\dexc(\pi)=k$, the permutation $\sigma$ having the desired 
parameters corresponds to the parallelogram polyomino obtained from the one associating to $\pi$ by reflection. 
Clearly, this operation does not change the diagonal lengths. Hence we have 
$|[\pi]|=|[\sigma]|$. Since $\exc$ and $\dexc$ are invariant on the classes 
this yields the proof.   
\end{bew}

The proof shows even more: not only the products of the numbers $a_1,\ldots,a_n$ 
counting possible letter exchanges are equally for $\pi$ and $\sigma$ but also 
the numbers themself, up to order. (If $a_i$ denotes the number of pairs 
$(i,\cdot)$ in $T_\pi$, and $a'_i$ denotes the number of pairs 
$(i,\cdot)$ in $T_\sigma$ then we have $a_i=a'_{n+1-i}$ for all $i$.) 
Consequently, we have
\bdpm
|\{\tau\in{\cal S}_n:\hat{\tau}=\pi,\;\inv(\tau)-\inv(\pi)=k\}|=
|\{\tau\in{\cal S}_n:\hat{\tau}=\sigma,\;\inv(\tau)-\inv(\sigma)=k\}|
\edpm 
for all $k$. (Recall that $\inv$ and $\dexc$ are identical statistics over ${\cal 
B}_n$.) We obtain the following result that was proved analytically in \cite[Cor. 
12]{clarke etal}.

\begin{cor}
$S_n^{(\exc,\inv)}(e,k)=S_n^{(\exc,\inv)}(n-1-e,n-1-2e+k)$ for all $e$ and $k$.
\end{cor}
\vspace*{0.75cm}

\centerline{\large\sc Acknowledgement}
\vspace*{0.5cm}

The most results in this paper are taken from my doctoral thesis \cite{reifegerste1}, prepared 
under the expert supervision of Christine Bessenrodt to whom I would like to 
express my gratitude. 
\vspace*{0.75cm}


\centerline{\large\sc References}
\vspace*{0.5cm}

\begin{enumbib}

\bibitem{barcucci etal1}
E. Barcucci, A. Del Lungo, E. Pergola, and R. Pinzani, 
From Motzkin to Catalan permutations, 
{\it Discrete Math.} {\bf 217} (2000), 33-49.

\bibitem{barcucci etal2}
E. Barcucci, A. Del Lungo, E. Pergola, and R. Pinzani, 
Some permutations with forbidden subsequences and their inversion number, 
{\it Discrete Math.} {\bf 234} (2001), 1-15.

\bibitem{biane}
P. Biane, 
Permutations suivant le type d'exc\'{e}dance et le nombre d'inversions, et 
interpr\'{e}tation combinatoire d'une fraction continue de Heine, 
{\it Europ. J. Comb.} {\bf 14} (1993), 277-284.

\bibitem{billey etal}
S. C. Billey, W. Jockusch, and R. P. Stanley, 
Some Combinatorial Properties of Schubert Polynomials, 
{\it J. Alg. Comb.} {\bf 2} (1993), 345-374.

\bibitem{clarke etal}
R. J. Clarke, E. Steingr\'{\i}msson, and J. Zeng, 
New Euler-Mahonian statistics on permutations and words, 
{\it Adv. Appl. Math.} {\bf 18} (1997), 237-270.

\bibitem{delest}
M. P. Delest, 
Polyominoes and animals: some recent results, 
{\it J. Math. Chem.} {\bf 8} (1991), 3-18.

\bibitem{delest-viennot}
M. P. Delest and X. G. Viennot, 
Algebraic languages and polyominoes enumeration, 
{\it Theoretical Comp. Sci.} {\bf 34} (1984), 169-206.

\bibitem{deutsch}
E. Deutsch, Problem 10816, 
{\it Amer. Math. Monthly} {\bf 107}, no. 7 (2000), 652.

\bibitem{deutsch-shapiro}
E. Deutsch and L. W. Shapiro,  
A survey of the Fine numbers,
{\it Discrete Math.} {\bf 241} (2001), 241-265.

\bibitem{donaghey-shapiro}
R. Donaghey and L. W. Shapiro, 
Motzkin numbers,
{\it J. Comb. Theory Ser. A} {\bf 23} (1977), 291-301.

\bibitem{foata}
D. Foata, 
Rearrangements of Words,
in M. Lothaire, {\it Combinatorics on Words}, Cambridge University Press, 1983.

\bibitem{foata-zeilberger}
D. Foata and D. Zeilberger, 
Denert's permutation statistic is indeed Euler-Mahonian,
{\it Stud. Appl. Math.} {\bf 83} (1990), 31-59.

\bibitem{francon-viennot}
J. Fran\c{c}on and X. G. Viennot, 
Permutations selon les pics, creux, doubles mont\'{e}es, doubles 
descentes, nombres d'Euler, nombres de Genocchi,
{\it Discrete Math.} {\bf 28} (1979), 21-35.

\bibitem{krattenthaler}
C. Krattenthaler, 
Permutations with restricted patterns and Dyck paths,
{\it Adv. Appl. Math.} {\bf 27} (2001), 510-530.

\bibitem{nazarov-tarasov}
M. Nazarov and V. Tarasov, 
On irreducibility of tensor products of Yangian modules associated with skew Young diagrams, 
{\it Duke Math. J.} {\bf 112}, (2002), 343-378.

\bibitem{reifegerste1}
A. Reifegerste, 
Differenzen in Permutationen: \"Uber den Zusammenhang von Permutationen, 
Polyominos und Motzkin-Pfaden,  
Ph.D. Thesis, University of Magdeburg, 2002.

\bibitem{reifegerste}
A. Reifegerste, 
On the diagram of $132$-avoiding permutations, 
preprint math.CO/0208006, 2002.

\bibitem{robertson etal}
A. Robertson, D. Saracino, and D. Zeilberger, 
Refined restricted permutations, 
preprint math.CO/0203033, 2002.

\bibitem{simion}
R. Simion, 
Combinatorial Statistics on Non-crossing Partitions, 
{\it J. Combin. Theory Ser. A} {\bf 66}, (1994), 270-301.

\bibitem{stanley1}
R. P. Stanley, 
{\it Enumerative Combinatorics Volume II},
Cambridge University Press, 1999.

\bibitem{stanley2}
R. P. Stanley, 
The Rank and Minimal Border Strip Decompositions of a Skew Partition, 
preprint, version of 29 August, 2002.

\bibitem{sulanke}
R. A. Sulanke, 
Counting Lattice Paths by Narayana Polynomials, 
{\it Electron. J. Comb.} {\bf 7} (2000), R40.

\bibitem{west}
J. West,  
Generating trees and the Catalan and Schr\"oder numbers,
{\it Discrete Math.} {\bf 146} (1995), 247-262.

\end{enumbib}
 
\end{document}